\documentclass[11pt]{amsart}

\usepackage{amssymb}
\usepackage{amsmath}

\begin{document}
\newtheorem{theo}{Theorem}[section]
\newtheorem{prop}[theo]{Proposition}
\newtheorem{lemme}[theo]{Lemma}
\newtheorem{cor}[theo]{Corollary}
\newtheorem{rem}[theo]{Remark}
\newtheorem{conj}[theo]{Conjecture}

\numberwithin{equation}{section}
\leftskip -2cm
\rightskip -2cm

\def\a{\alpha}
\def\ind{{\rm ind}}
\def\aa{{\mathbb A}}
\def\bc{\backslash}
\def\lg{L^2(GL_n(\F)Z(\aa_\F)\bc GL_n(\aa_\F),\o)}
\def\o{\omega}
\def\lgg{L^2(G'(F)Z(\aa)\bc G'(\aa);\o)}
\def\lra{\leftrightarrow}
\def\ccc{{\bf C}}
\def\jlr{{\bf JL}_r}
\def\ljr{{\bf LJ}_r}

\def\lge{L^2(GL_n(\E)Z(\aa_\E)\bc GL_n(\aa_\E),\o)}

\def\S{{\mathfrak S}}

\numberwithin{equation}{section}

\def\e{\varepsilon}
\def\f{{\mathcal F}}
\def\i{{\bf i}}
\def\k{\{1,2,...,k\}}

\def\n{\mathbb N}
\def\r{{\mathbb R}}
\def\s{\sigma}
\def\z{\mathbb Z}
\def\q{{\mathbb Q}}

\def\t{\times}
\def\tt{\times}
\def\cc{{\mathbb C}}

\def\lra{\leftrightarrow}
\def\bc{\backslash}
\def\cusp{{\mathcal C}}\def\tr{{\rm tr}}

\def\ge{G_E}
\def\gf{G_F}

\def\gaae{G(\aa)_E}
\def\gaaf{G(\aa)_F}

\def\nn{{\mathcal N}}

\def\nio{normalized $\s$-intertwining operator}

\def\d{\delta}
\def\D{\Delta}

\def\u{{\mathcal U}}
\def\g{\gamma}
\def\G{\Gamma}

\def\F{{\bf F}}
\def\E{{\bf E}}

\centerline{\huge Shintani relation for base change:}
\centerline{\huge unitary and elliptic representations}
\ \\

\centerline{A. I. Badulescu, G. Henniart}
\ \\

\centerline{\it To Jim Cogdell for his 60th birthday}
\ \\

{\small
\tableofcontents}

\ \\
{\bf Abstract.} Let $E/F$ be a cyclic extension of $p$-adic fields and $n$ a positive integer. Arthur and Clozel constructed a base change process $\pi\mapsto \pi_E$ which associates to a smooth irreducible representation of $GL_n(F)$ a smooth irreducible representation of $GL_n(E)$, invariant under $Gal(E/F)$. When $\pi$ is tempered, $\pi_E$ is tempered and is characterized by an identity (the Shintani character relation) relating the character of $\pi$ to the character of $\pi_E$ twisted by the action of $Gal(E/F)$. In this paper we show that the Shintani relation also holds when $\pi$ is unitary or elliptic. We prove similar results for the extension $\cc/\r$. As a corollary we show that for a cyclic extension $\E/\F$ of number fields
the base change for automorphic residual representations of $GL_n(\aa_\F)$  respects the Shintani relation at each place of $\F$.

\section{Introduction} Our goal in this paper is to extend results of Arthur and Clozel ([AC]) on local base change for tempered representations and global base change for cuspidal representations. 
In particular, for $p$-adic fields, we show that the Shintani character relation expresses base change for unitary or elliptic representations. 
Globally we construct base change for residual representations compatible with the local Shintani character relations.\\

Let us give some detail, first in the $p$-adic case, for a prime number $p$. 
Let $E/F$ be a cyclic extension of $p$-adic fields of degree $l$, $\s$ a generator of $Gal(E/F)$ and $n$ a positive integer. 
To each smooth irreducible representation $\pi$ of $GL_n(F)$, Arthur and Clozel attach a pair $(\pi_E,I_{\pi_E})$, called here {\bf base change of $\pi$}, where $\pi_E$ is a smooth irreducible representation of $GL_n(E)$ and $I_{\pi_E}$ is an isomorphism of $\pi_E$ onto $\pi_E^{\s}$. When $\pi$ is tempered, $\pi_E$ is tempered and is characterized by the Shintani character relation:
\begin{equation}\label{Shintani}
\chi_{\pi_E,\s}(g)=\chi_\pi(\nn g)
\end{equation}
for any $\s$-regular element $g$ of $GL_n(E)$: in this relation $\chi_\pi$ is the character of $\pi$, $\chi_{\pi_E,\s}$ is the twisted character of $\pi_E$ associated to the choice of $\s$, and $\nn g$ is an element of $GL_n(F)$, called the norm of $g$, conjugate in $GL_n(E)$ to $gg^\s g^{\s^2}...g^{\s^{l-1}}$. 
This result extends easily to essentially tempered representations and to Levi subgroups of $GL_n(F)$ instead of $GL_n(F)$.
In general, $\pi$ is the Langlands quotient of a parabolically induced representation from an essentially tempered representation $\tau$, of some Levi subgroup, and $\pi_E$ is then defined to be the Langlands quotient of the induced representation from $\tau_E$. 
Easy examples show that the Shintani relation does not hold for general $\pi$ and its base change (see {\bf Example} after the prop. \ref{liftprod}).
If the Shintani relation holds, we say that $\pi_E$ is a {\bf Shintani lift} of $\pi$.
Our main local result is:\\
\ \\
{\bf Theorem.} {\it When $\pi$ is unitary or elliptic, $\pi_E$ is a Shintani lift of $\pi$}.\\

We also prove the same result in the case of the Archimedean extension $\cc/\r$.\\

In the sequel, we say that $\pi$ has a Shintani lift to mean that $\pi$ and its base change $\pi_E$ verify the Shintani relation -- thus it is automatic when $\pi$ is tempered. We call Shintani lift or simply lift the process of showing that certain classes of representations of $GL_n(F)$ have a Shintani lift.

Our global result ({\bf Theorem E}, section \ref{resultslocal}) concerns residual automorphic representations. We consider a cyclic extension $\E/\F$ of degree $l$ of number fields, and an automorphic discrete series representation $\pi$ of $GL_n(\aa_\F)$. If $v$ is a place of $\F$ and $w$ a place of $\E$ above $v$, the component $\pi_v$ of $\pi$ at $v$ has a Shintani lift $\pi'_w$ to $GL_n(\E_w)$ and we form the admissible representation $\pi':=\otimes_w \pi'_w$ of $GL_n(\aa_\E)$. When $\pi$ is cuspidal, Arthur and Clozel showed that $\pi'$ is automorphic, induced from cuspidal. We extend it to the residual case and show that if $\pi$ is residual, then $\pi'$ is automorphic, and $\pi'$ is parabolically irreducibly induced from a tensor product of residual representations. 

Let us comment on the proofs which are gathered in section \ref{proofs}.

The Theorem for elliptic representations is deduced by local methods from the case of square integrable representations (done in [AC]) using the explicit description of elliptic representations (see section \ref{sectelliptic}). 

To treat unitary representations, an important local step is the following, proved in section \ref{proofs}:\\
\ \\
{\bf Proposition.} {\it Let $\pi_i$ be smooth irreducible representations of $GL_{n_i}(F)$, $i=1,2$. Assume that $\pi_i$ has Shintani lift $\pi_{i,E}$ and that the parabolically induced representations $\pi_1\tt \pi_2$ and $\pi_{1,E}\tt \pi_{2,E}$ are irreducible. 
Then the base change of $\pi_1\tt \pi_2$ is $\pi_{1,E}\tt \pi_{2,E}$. The representation $\pi_{1,E}\tt \pi_{2,E}$ is also a Shintani lift for $\pi_1\tt \pi_2$.}\\

The notation $\pi_1\t\pi_2$ stands for the representation of $GL_{n_1+n_2}(F)$ obtained by normalized parabolic induction from $\pi_1\otimes \pi_2$,  with respect to the upper triangular parabolic subgroup with Levi subgroup the group of block diagonal matrices of size $n_1$ and $n_2$. This Proposition implies a short proof of the Theorem when $\pi$ is unitary and spherical, a result proved in [AC] by global means. It also implies that to prove the Theorem for general unitary irreducible representations, it is enough to prove it for the so-called Speh representations $u(\delta,k)$ (see section \ref{classification} for the definition). To prove the Theorem for Speh representations, we follow the local-global method of Arthur and Clozel and use the trace formula ([AC] Chapter 2). For such a Speh representation $\pi=u(\delta,k)$ there is :

(i) a cyclic extension $\E/\F$ of number fields giving the extension $E/F$ at a place $v$ of $\F$, 

(ii) an automorphic discrete series representation $\Pi$ of $GL_n(\aa_\F)$ such that $\Pi_v\simeq \pi$ and such that at all places $v'$ of $\F$ different from $v$,  $\Pi_{v'}$ has a Shintani lift.\\
From the trace formula it follows that $\pi$ has a Shintani lift. The same type of proof works for the Archimedean case. Once the local lift is proved for all unitary representations, the global lift for residual automorphic representations is obtained by the trace formula of [AC] and arguments of compatibility between the local and global settings. 

We want to thank the referee warmly for his careful reading and his many remarks which helped us to substantially improve the quality of this presentation.

\section{Notation and basic facts (local)}\label{notation}

In the following sections of chapter \ref{notation}, except the last one \ref{archim}, $F$ will be a $p$-adic field. We write $O_F$ for the ring of integers of $F$ and $q_F$ for the cardinality of the residue field. We fix a
uniformizer $\pi_F$ of $F$ and let $|\ \ |_F$ be the absolute value
on $F$ defined by $|\pi_F|_F=q_F^{-1}$. We consider complex smooth representations of linear groups over $F$, which we simply call ``representations".

\subsection{Classifications}\label{classification}
Let $n$ be a positive integer. Put  $G_F=GL_n(F)$. Let $Z_F$ be the center of $G_F$. Let $K_F:=GL_n(O_F)$ and fix Haar measures on $G_F$, resp. $Z_F$, such that $vol(K_F)=1$, resp. $vol(Z_F\cap K_F)=1$. Let $B_F$ be the Borel subgroup made of upper triangular matrices and $U_F$ the unipotent subgroup of $B_F$ made of upper triangular matrices with $1$ on the diagonal.\\

For any $k\in\n^*$, let $\nu$ denote the character of $GL_k(F)$ given by the composition of the norm $|\ \ |_F$ with the determinant map.
The twist of a representation $\pi$ (of $G_F$) with a character $\chi$ (of $G_F$) will be written simply $\chi\pi$ (instead of $\chi\otimes \pi$). 
Moreover, {\it when we write $\chi\pi$ for a character $\chi$ of $F^\t$ and a representation $\pi$ of $G_F$ we mean $(\chi\circ\det)\pi$}.\\

A smooth representation $\rho$ of $G_F$ is {\bf cuspidal} if $\rho$ is irreducible and
has a non-zero coefficient with compact support modulo
$Z_F$. A smooth representation $\d$ of $G_F$ is {\bf
square integrable} if $\d$ is irreducible, unitary, and has a non-zero coefficient
which is square integrable over $G_F/Z_F$. 
An {\bf essentially square integrable} representation is the twist of a square integrable 
representation with a character.
Any essentially square integrable representation $\d$ may be written as $\d=\nu^{e(\d)}\d^u$ where $e(\d)$ is a real number, and $\d^u$ is a square integrable  representation. 
Then $\d$ determines $e(\d)$ -- which is called the exponent of $\d$ -- and $\d^u$.\\

A {\bf standard Levi subgroup} of $G_F$ is a subgroup $L$ of block diagonal matrices of given sizes. If $n_1,n_2,...,n_k$ (where $\sum_{i=1}^k n_i=n$) are the sizes of blocks,  then $L$ is identified with the product $\prod_{i=1}^k GL_{n_i}(F)$ (which is $GL_{n_1}(F)\t GL_{n_2}(F)\t ...\t GL_{n_k}(F)$ in this order). 
We denote by $P_L$ 
the parabolic subgroup generated by $L$ and $B_F$. The
definitions of cuspidal representation and square integrable representation extend to $L$ in the obvious way. 
Here $\ind_L^{G_F}$ will denote the normalized parabolic induction from 
$(L,P_L)$ to $G_F$. Then, if $\pi_i$ is an admissible representation of
$G_{n_i}(F)$ for $1\leq i\leq k$ we write
$\pi_1\times\pi_2\times...\times \pi_k$ for
$\ind_L^{G_F} \pi_1\otimes\pi_2\otimes...\otimes \pi_k$.  We call $\pi_1\times\pi_2\times...\times \pi_k$ the product of the representations $\pi_i$ and we sometimes denote it $\prod_{i=1}^k\pi_i$, not forgetting that the product is not commutative.\\

Let us recall the Bernstein-Zelevinsky classification of essentially square integrable representations ([BZ], [Ze]). If $\d$ is a square integrable representation of $GL_n(F)$, there exists a pair $(k,\rho)$, where $k$ is a divisor of $n$ and $\rho$ is a unitary cuspidal representation of $GL_\frac{n}{k}(F)$ such that $\d$ is isomorphic to the unique irreducible subrepresentation $Z(\rho,k)$ of $\nu^{\frac{k-1}{2}}\rho\t \nu^{\frac{k-1}{2}-1}\rho \t ...\t \nu^{-\frac{k-1}{2}}\rho$ (this induced representation also has a unique irreducible quotient, its Langlands quotient, defined below). The integer $k$ and the isomorphism class of $\rho$ are determined by the isomorphism class of $\delta$.

Let $\delta$ be an {\it essentially} square integrable representation. Then there exist a positive integer $k$, $k|n$, and a cuspidal representation $\rho$ of $GL_{\frac{n}{k}}(F)$ such that $\delta$ is  the unique irreducible subrepresentation of 
$\nu^{k-1}\rho\times \nu^{k-2}\rho\times ...\times\rho$. The set 
$\{\rho,\nu \rho, ... ,\nu^{k-1}\rho\}$ is a {\bf Zelevinsky segment}, the Zelevinsky segment of $\d$. The integer $k$ is its {\bf length}.\\

If $\d$ is an irreducible unitary representation of a standard Levi subgroup $L$, in particular if $\d$ is square integrable, then
 the induced representation $\ind_L^{G_F}\d$ is irreducible ([Be]). 
 We define here a {\bf tempered} representation to be an (irreducible) representation of the form $\ind_L^{G_F}\d$, 
 where $\delta$ is a square integrable representation of $L$.
An {\bf essentially tempered} representation is the twist of a tempered representation with a character.
Then any essentially tempered representation $\tau$ may be written as $\tau=\nu^{e(\tau)}\tau^u$ where $e(\tau)$ is a real number, and $\tau^u$ is a tempered representation. 
Then $\tau$ determines $e(\tau)$ -- which is called the exponent of $\tau$ -- and $\tau^u$.\\

We recall the Langlands classification ([Re] VII.4.2).
Let $\tau_1,\tau_2,...,\tau_k$ be tempered representations of groups $GL_{n_i}(F)$ and $\a_1>\a_2>...>\a_k$ real numbers. 
We say that the essentially tempered representations $\nu^{\a_i}\tau_i$ are in standard or strictly decreasing order. 
Then $\prod_{i=1}^k\nu^{\a_i}\tau_i$ has a unique irreducible quotient, called the {\bf Langlands quotient} and 
denoted here  $Lg(\nu^{\a_1}\tau_1,\nu^{\a_2}\tau_2,...,\nu^{\a_k}\tau_k)$.
Every irreducible representation $\pi$ of $G_n(F)$ is isomorphic with such a $Lg(\nu^{\a_1}\tau_1,\nu^{\a_2}\tau_2,...,\nu^{\a_k}\tau_k)$ such that $k$, $\a_1,\a_2,...,\a_k$ and the isomorphism classes of $\tau_1,\tau_2,...,\tau_k$ are  determined by the isomorphism class of $\pi$. 
We may extend this definition to the case when $\a_1\geq\a_2\geq...\geq\a_k$, and this will be convenient, for example, for the proof of Proposition \ref{nioprod}. We say then that the essentially tempered representations $\nu^{\a_i}\tau_i$ are in decreasing order.
Set
$\a_{i_1}=\a_1$ and let $\a_{i_2} > \a_{i_3} > ... > \a_{i_r}$ be the subsequence of $\a_1,\a_2,...,\a_k$ made of all elements $\a_i$ such that $\a_{i-1}> \a_i$.
Set $\tau'_{i_j}:=\prod_{i_j\leq i\leq i_{j+1}-1}\tau_i$ ($i_{r+1}:=k+1$). Then $\tau'_{i_j}$ is an irreducible tempered representation and we have $\prod_{i=1}^k\nu^{\a_i}\tau_i \simeq \prod_{j=1}^r\nu^{\a_{i_j}}\tau'_{i_j}$. We set by definition  $Lg(\nu^{\a_1}\tau_1,\nu^{\a_2}\tau_2,...,\nu^{\a_k}\tau_k) := Lg(\nu^{\a_{i_1}}\tau'_{i_1},\nu^{\a_{i_2}}\tau'_{i_2},...,\nu^{\a_{i_r}}\tau'_{i_r})$.
Moreover, if an irreducible representation $\pi$ of $G_F$ is given, then $\pi$ is isomorphic to some $Lg(\nu^{\a_1}\tau_1,\nu^{\a_2}\tau_2,...,\nu^{\a_k}\tau_k)$ with $\tau_i$ square integrable and $\a_1\geq \a_2\geq ...\geq \a_k$. If $\tau_i$ are given essentially tempered representations, we define $Lg(\tau_1,\tau_2,...,\tau_k)$ or $Lg(\tau_1\t\tau_2\t...\t\tau_k)$
as being the Langlands quotient when we order $\tau_i$ in such an order that their exponents are in decreasing order, and this is independent of the order we choose with that property.

If $\tau$ is a tempered representation and
$k\in\n^*$, we will write $u(\tau,k)$ for $Lg(\nu^{\frac{k-1}{2}}\tau,
\nu^{\frac{k-1}{2}-1}\tau,...,\nu^{-\frac{k-1}{2}}\tau)$ (when $\tau$ is square integrable, these are Speh representations). Moreover, if
$\a\in]0,\frac{1}{2}[$, we let $\pi(u(\tau,k),\a)$ be the
representation $\nu^{\a}u(\tau,k)\times \nu^{-\a}u(\tau,k)$. It is an irreducible representation.\\

The irreducible unitary representations of the groups $GL_n(F)$ have been
classified by Tadi\'c ([Ta1]), using also the already quoted result of 
Bernstein ([Be]). We describe Tadi\'c's result. Let $\u$ be the set of isomorphism classes of representations
$u(\d,k)$, $\pi(u(\d,k),\a)$ where $\d$ is a square integrable representation, $k$
is a positive integer and $\a\in]0,\frac{1}{2}[$. Then all the
representations in $\u$ are irreducible and unitary. Any product of
representations in $\u$ is irreducible. Any irreducible unitary
representation $\pi$ of any $GL_n(F)$, $n\in\n^*$, is isomorphic to such a product
of representations from $\u$, and the factors of the product are
determined by $\pi$ up to isomorphism and permutation.\\

Let  $\psi$ be a non-trivial additive character of $F$; we define $\Theta_\psi:U_F\to \cc^\t$ by $\Theta_\psi(u)=\psi(\sum_{i=1}^{n-1} u_{i,i+1})$ if $u=(u_{i,j})_{1\leq i,j\leq n}$. A {\bf Whittaker functional} on a representation $(\pi,V)$ of $GL_n(F)$ is a linear map $\lambda:V\to\cc$ such that $\lambda(\pi(u)v)=\Theta_\psi(u)\lambda(v)$ for $v\in V$, $u\in U_F$. If $\pi$ is irreducible, the space of Whittaker functionals $W(\pi,\psi)$ is of dimension $0$ or $1$. We say that a smooth representation $\pi$ is {\bf generic} if $\pi$ is irreducible and the dimension of $W(\pi,\psi)$
is $1$, a condition which does not depend of the choice of $\psi$. (For $n=1$, $U_F$ is trivial, so all the characters of $GL_1(F)$ are generic.)

According to [Ze] Theorem 9.7, an irreducible representation $\gamma$ is generic if and only if $\gamma$ is isomorphic to an irreducible product
$\prod\nu^{\a_i}\d_i$ with $\a_i\in\r$ and $\d_i$ are square integrable representations (in [Ze], generic representations are called non-degenerate). 
If a generic representation $\gamma$ is the Langlands quotient of a representation induced from an essentially tempered representation $t$ of a standard Levi subgroup, it follows that $\gamma$ is the full induced representation from $t$. 
If, moreover, $\gamma$ is unitary, Tadi\'c's classification implies that $\gamma$ is an irreducible product 
$$\d_1\t\d_2\t...\t\d_k  \t \pi(\d'_1,\a_1)\t \pi(\d'_2,\a_2)\t...\t \pi(\d'_l,\a_l),$$ 
where $\d_1,\d_2,...,\d_k, \d'_1,\d'_2,...,\d'_l$ are square integrable representations and $\a_1,\a_2,...,\a_l\in ]0,\frac{1}{2}[$. Note that such a product is always irreducible, unitary and generic. Hence a product of unitary generic representations is a unitary generic representation. Note also that tempered representations are generic and unitary.

If $\gamma$ is a unitary generic representation as before and
$k\in\n^*$, then the induced representation $\nu^{\frac{k-1}{2}}\g\t\nu^{\frac{k-1}{2}-1}\g\t...\t\nu^{-\frac{k-1}{2}}\g$ has a unique
irreducible quotient which we denote $u(\g,k)$ (we recall that $u(\g,k)$ has already been defined when $\g$ is tempered). With the notation at the beginning of the last paragraph, we have then
$u(\g,k)\simeq \prod\nu^{\a_i}u(\d_i,k)$ (see [Ba2] 4.1, for example). In
particular, Tadi\'c's classification implies that $u(\g,k)$ is unitary. Notice that the local component of a global cuspidal representation is unitary generic, so the local component of an automorphic residual representation, according to the Moeglin and Waldspurger classification ([MW1], recalled below, section \ref{notglobal}), is always of type $u(\gamma,k)$ for some unitary generic representation $\gamma$.\\

A representation $\pi$ of $G_F$ is said to be {\bf spherical} if $\pi$ is irreducible and has a non-zero vector fixed under $K_F=GL_n(O_F)$.
Let $\pi=Lg(\d_1,\d_2,...,\d_k)$ be an irreducible representation of $G_F$, where $\d_1,\d_2,...,\d_k$ are essentially square integrable representations. Then $\pi$ is spherical if and only if all the $\d_i$, $1\leq i\leq k$, are unramified characters of $GL_1$.

\subsection{The \nio}\label{niosection}

Let $E$ be a cyclic extension of $F$ and $l:=[E:F]$. Let $N_{E/F}:E\to F$ be the norm of the extension. Fix a generator $\s$ of $Gal(E/F)$. Put $G_E=GL_n(E)$.  The Haar measure on $G_E$ is defined in the same way as on $G_F$. If $s\in Gal(E/F)$ and $x\in G_E$, then we write $s(x)$ or $x^s$ for the matrix obtained from $x$ by applying $s$ to all the coefficients. Then $x^{ss'}=(x^{s'})^s$, $\det(x^s)=s(\det(x))$ and  $x\mapsto x^s$ is a group automorphism of $G_E$.

If $(\Pi,V)$ is a representation of $G_E$ and $s\in Gal(E/F)$, we write $(\Pi^s,V)$ for the representation of $G_E$ given by $g\mapsto \Pi(g^{s})$. We then say $\Pi$ is $s$-{\bf stable} if  $\Pi^s\simeq \Pi$. Let us say a few words about the behavior of this Galois action on representations with respect to the induction functor. Let $(\Pi,V)$ be a representation of $G_E$ by right translations in some space of functions on $G_E$ and $s\in Gal(E/F)$. For $f\in V$ define $f_s$ by $x\mapsto f(s(x))$. Let $sV$ be the space $\{f_s,\ f\in V\}$ and let $(s\Pi,sV)$ be the representation of $G_E$ by right translations in $sV$. By abuse, we denote sometimes $s:V\to sV$ the map $f\mapsto f_s$, a notation used also in [AC]. 
Assume $\pi=\ind_L^{G_E}\tau$ is an induced representation, where $(\tau,W)$ is a representation of a standard Levi subgroup $L$ of  $G_E$. Then $\ind_L^{G_E}(\tau^s) =s\pi$. Indeed, by definition, $\pi$ is the representation by right translations in the space 
  $$V:=\{f:G_E\to W{\text{ smooth }}, f(pg)=\d^{\frac{1}{2}}(p)\tau (p)f(g) \forall g\in G_E, p\in P_L\}$$
and $\ind_L^{G_E}\tau^s$ is the representation of $G_E$ by right translations in the space:
  $$U:=\{f:G_E\to W {\text{ smooth }}, f(pg)=\d^{\frac{1}{2}}(p)\tau (s(p))f(g) \forall g\in G_E, p\in P_L\},$$
where $\d$ is the modulus character of the parabolic subgroup $P_L$ from which we induce. Notice that the space of $\tau^s$ is the same space $W$ as for $\tau$, but $V\neq U$ (unless, for example, $\tau^s=\tau$).
Now $f\mapsto f_s$ is a bijection from $V$ to $U$ which intertwines $\pi^s$ with $s\pi$, hence an identification $U=s V$ and $\ind_L^{G_E}(\tau^s)=s\pi$. 

Let 
$I :W\to W$ be an intertwining operator between $\tau$ and $\tau^s$. Parabolic induction is a functor, and we let
 $\ind_L^{G_E}I:V\to U$ be the intertwining operator (between $\pi$ and $s\pi$) parabolically induced from $I$. 
We let $I_s(\pi):=s^{-1}\circ \ind_L^{G_E}I$, where $s:f\mapsto f_s$.\\
\ \\
{\bf Definition.} We say  $I_s(\pi)$ {\bf is the $s$-operator on $\pi$ obtained from $I$ by the parabolic induction procedure}.\\

The operator $I_s(\pi)$ is an operator on the space of $\pi$ which intertwines $\pi$ and $\pi^s$.
If, for example, $n=2$, $L=GL_1(E)\t GL_1(E)$ 
and $\tau$ is an  unramified character of $L$, if $I$ is trivial, then $V=U$, but $I_s(\tau)$ is equal to 
$f(^.)\mapsto f(s^{-1}(^.))$, which is not identity (unless $E=F$).\\
\ \\
{\bf Important.} Most often we will meet the following situation: $\pi=\ind_L^{G_E}\tau$ has an irreducible subquotient $\pi_0$ of our interest, which is multiplicity one and $s$-stable. Then the $s$-operator $I_s(\pi)$ on $\pi$ obtained from $I$ by the parabolic induction procedure induces by the multiplicity one property (see section {\it Group with automorphism} of the Appendix) an operator $I_s(\pi_0)$ on $\pi_0$, well defined in the sense that it does not depend on the way we realize $\pi_0$ as a subquotient of $\pi$.\\
\ \\
{\bf Definition.} We say $I_s(\pi_0)$ is {\bf the $s$-operator on $\pi_0$ obtained from $I$ by p.i.m.o.} (p.i.m.o. stands for parabolic induction and multiplicity one).\\

Let $L'$ be a standard Levi subgroup of $G_E$ such that $L\subset L'$. Let
$$V':=\{f:L'\to W{\text{ smooth}}, f(pg)=\d^{\frac{1}{2}}(p)\tau (p)f(g) \forall g\in L', p\in P_{L}\cap L'\}.$$
The representation $\tau':=\ind_L^{L'}\tau$ is the representation by right translations of $L'$ in $V'$. 
Then the representation $\ind_{L'}^{G_E}\tau'$ is the representation by right translations of $G_E$ in the space
$$V'':=\{f:G_E\to V'{\text{ smooth}}, f(pg)=\d^{\frac{1}{2}}(p)\tau' (p)f(g) \forall g\in G_E, p\in P_{L'}\}.$$
It is known, and easy to check, that the map $h:V''\to V$, defined by $h(f)=(g\mapsto f(g)(1))$ for $f\in V''$ and $g\in G_E$, is an isomorphism of representations from  $\ind_{L'}^{G_E}\tau'$ to $\pi=\ind_{L}^{G_E}\tau$ (this is the transitivity of the parabolic induction functor).

With these notation, if $\tau$ is $s$-stable and $I$ intertwines $\tau$ and $\tau^s$, we have the transitivity property:

\begin{prop}\label{transitivite} 
(a) One has $h\circ (I_s(\tau'))_s(\pi)=I_s(\pi)\circ h$.

(b) Let $\pi_0$ be an irreducible $s$-stable subquotient of $\pi$ of multiplicity one.
Let $\tau'_0$ be the irreducible subquotient of $\tau'$ such that $\pi_0$ is a subquotient of $\ind_{L'}^{G_E}\tau'_0$. If $\tau'_0$ is $s$-stable we have:
$$I_s(\pi_0)=(I_s(\tau'_0))_s(\pi_0).$$
\end{prop}
\ \\
Claim (b) will be important for crucial proofs in this article.\\
\ \\
{\bf Proof.} (a) follows from simple verification starting with the definition of $h$.

(b) 
Let 
$$0\longrightarrow W\longrightarrow U\to \tau'_0\longrightarrow 0$$ 
be an exact sequence of representations, where $(U,W)$ is the maximal pair of subrepresentations of $\tau'$ such that $U/W\simeq \tau'_0$. As explained in the Appendix, $U$ and $W$ are stable by $I_s(\tau')$. The parabolic induction functor is exact ([Re], Prop. II.2.2) and we obtain:
$$0\ \ \ \ \rightarrow\ \ \ \  \ind_{L'}^{G_E} W\ \ \ \ \rightarrow\ \ \ \  \ind_{L'}^{G_E} U\ \ \ \ \xrightarrow{F}\ \ \ \  \ind_{L'}^{G_E} \tau'_0\ \ \ \ 
\rightarrow 0.$$

Now $\pi_0$ is a subquotient of $\ind_{L'}^{G_E}\tau'_0$ of multiplicity one, and we let $(u,w)$ be the maximal pair of subrepresentations of $\ind_{L'}^{G_E}\tau'_0$ such that $u/w\simeq \pi_0$ (see the Appendix). We end up with a chain of inclusions 
$$\ind_{L'}^{G_E} W\ \ \subset\ \  F^{-1}(w)\ \ \subset\ \  F^{-1}(u)\ \ \subset\ \   \ind_{L'}^{G_E} U,$$
such that the isomorphism 
$$\ind_{L'}^{G_E} U/ \ind_{L'}^{G_E} W\simeq \ind_{L'}^{G_E}\tau'_0$$ 
sends 
$$F^{-1}(u)/\ind_{L'}^{G_E} W \text{ onto } u$$ 
and 
$$F^{-1}(w)/\ind_{L'}^{G_E} W \text{ onto } w.$$

By (a), the operator on $\ind_{L'}^{G_E} U$ obtained by restriction from $I_s(\pi)$ equals $J_s(\ind_{L'}^{G_E} U)$ where $J$ is obtained on $U$ by restriction from $I_s(\tau')$. Now $F^{-1}(u)$ is the maximal submodule of $\ind_{L'}^{G_E} U$ admitting $\pi_0$ as a quotient (indeed, the image by $F$ of such a submodule is included in $u$ by Proposition \ref{subquot} (a)).
The result then follows from the section {\it Group with automorphism} of the Appendix.\qed
\ \\

Let $\Pi$ be an irreducible $\s$-stable representation of $G_E$. Following [AC], we want to produce a {\it canonical} isomorphism $I=I_\Pi$ of $\Pi$ onto $\Pi^\s$; in particular, if $\phi$ is an isomorphism of $\Pi$ onto another representation $\Pi'$ of $G_\E$ we shall have by construction $I_{\Pi'}\circ\phi=\phi\circ I_\Pi$.

Let us first treat the case where $\Pi$ is {\it generic}. Choosing a non-trivial additive character $\psi$ of $F$, we put $\psi_E:=\psi\circ\tr_{E/F}$. 
Because $\psi_E$ is invariant under $\s$, we have $W(\Pi^\s,\psi_E)=W(\Pi,\psi_E)$ and any isomorphism of $\Pi$ onto $\Pi^\s$ induces an automorphism (actually, a non-zero homothety) of $W(\Pi,\psi_E)$. Consequently, there is a unique isomorphism $I_\Pi^{gen}$ of $\Pi$ onto $\Pi^\s$ inducing identity on $W(\Pi,\psi_E)$. This is the {\bf \nio} for generic representations.

Let $a\in F^\t$ and put $\psi^a:x\mapsto \psi(ax)$; then the map $\lambda\mapsto \lambda\circ \Pi({\rm diag}(a^{n-1}, a^{n-2},...,a,1))$ gives an isomorphism of $W(\Pi,\psi_E)$ onto $W(\Pi,\psi_E^a)$; since $a$ is fixed by $\s$, we see that $I_\Pi^{gen}$ does not depend on the choice of $\psi$. 
If $\phi$ is an isomorphism of $\Pi$ onto $\Pi'$ then $\lambda'\mapsto \lambda'\circ\phi$ is an isomorphism of $W(\Pi',\psi_E)$ onto $W(\Pi,\psi_E)$.
So {\it the operators $I_\Pi^{gen}$ are compatible with isomorphisms}.\\

Let now $\Pi$ be general, irreducible and $\s$-stable. Then $\Pi$ is in some isomorphism class $Lg(\Pi_1,\Pi_2,...,\Pi_k)$, where $\Pi_1,\Pi_2,...,\Pi_k$ are essentially tempered - hence generic - in standard order $e(\Pi_1)>e(\Pi_2)>...>e(\Pi_k)$.  As $\Pi^\s$ is then in the class $Lg(\Pi_1^\s,\Pi_2^\s,...,\Pi_k^\s)$, by unicity we get that $\Pi_i$ is $\s$-stable. Then the {\bf \nio\ $I_\Pi$ of $\Pi$} is by definition the $\s$-operator on $\Pi$ obtained from $I_{\Pi_1}^{gen}\otimes I_{\Pi_2}^{gen}\otimes ...\otimes I_{\Pi_k}^{gen}$ and $\s$ by p.i.m.o.. We have the:

\begin{lemme}\label{nio}
Let $\Pi$, $\Pi'$ be irreducible $\s$-stable isomorphic representations of $GL_n(E)$. If $h:\Pi\to \Pi'$ is an isomorphism, then 
$I_{\Pi'}=h I_\Pi h^{-1}$.
\end{lemme}
\ \\
{\bf Proof.} If $\Pi=Lg(\Pi_1,\Pi_2,...,\Pi_k)$, we fix a surjective intertwining operator $f:\Pi_1\t\Pi_2\t ...\t\Pi_k\to \Pi$ and we obtain 
a surjective intertwining operator $h\circ f:\Pi_1\t\Pi_2\t ...\t\Pi_k\to \Pi'$. The Lemma follows from the fact that $I_{\Pi}$ and $I_{\Pi'}$ do not depend on the choice of $f$.\qed

\def\l{\lambda}
\begin{lemme}\label{gengen}
Let $\Pi$ be a generic $\s$-stable representation of $G_E$. Then $I_\Pi^{gen}=I_\Pi$.
\end{lemme}

\ \\
{\bf Proof.} That is Lemma 2.1 in [AC] I which, however, offers only a hint of the proof. One can reason as follows. Write $\Pi$ as a Langlands quotient $Lg(\Pi_1,\Pi_2,...,\Pi_k)$ as above. If $\Pi$ is generic, in fact, the product $\Pi_1\t \Pi_2\t ...\t\Pi_k$ is irreducible and $Lg(\Pi_1,\Pi_2,...,\Pi_k)=\Pi_1\t \Pi_2\t ...\t\Pi_k$. We have noticed, just after having defined $I_\Pi^{gen}$, that these operators are compatible with isomorphisms, and so is the operator $I_\Pi$, by Lemma \ref{nio}. So we may assume $\Pi=\Pi_1\t \Pi_2\t ...\t\Pi_k$.

Choose non-zero Whittaker functionals $\l_i\in W(\Pi_i,\psi_E)$ for $i=1,2,...,k$. By [JS1] chapter 3, we have on $\Pi_1\t \Pi_2\t ...\t\Pi_k$ a Whittaker functional $\Lambda$ given by ([JS1], formula (2) chapter 3):
$$\Lambda(f)=\int_{U_E}\l(f(u))\overline{\Theta_{\psi_E}(u)}du,$$
where $\l=\l_1\otimes \l_2\otimes...\otimes \l_k$, $f$ is a function in the space of $\Pi_1\t\Pi_2\t ...\t\Pi_k$ (in particular, a map from $G_E$ to $V_1\otimes V_2\otimes ...\otimes V_k$, where $V_i$ is the space of $\Pi_i$), 
 and  $du$ is a Haar measure on $U_E$; by {\it loc.cit.} the integral is always convergent.

Now, for $u\in U_E$, we have:
$$\l(I_\Pi(f)(u))=\l(I(f(\s^{-1}(u))))$$
where $I$ is the product $I_{\Pi_1}^{gen}\t I_{\Pi_2}^{gen}\t ...\t I_{\Pi_k}^{gen}$, so
$$\l(I_\Pi(f)(u))=\l(f(\s^{-1}(u)))$$
by the definition of $I_{\Pi_i}^{gen}$. We then get $\Lambda\circ I_\Pi=\Lambda$ from the $\s$-invariance of $\Theta_{\psi_E}$ and $du$. The result follows.\qed
\ \\
\begin{prop}\label{transitivity}
Let $L$ be a standard Levi subgroup of $G_E=GL_n(E)$, $\gamma$ a $\s$-stable generic representation of $L$ and $I_\gamma=I^{gen}_\gamma$ the \nio\ of $\gamma$. 
Then $\ind_L^{G_E} \gamma$ has a unique generic irreducible subquotient $\Pi$. Moreover, $\Pi$ is $\s$-stable and
if $I_{\gamma,\s}(\Pi)$ is the $\s$-operator on $\Pi$ obtained by p.i.m.o. then $I_{\gamma,\s}(\Pi)=I_{\Pi}^{gen}$.
\end{prop}
\ \\
{\bf Proof.} The induced representation $\ind_L^{G_E} \gamma$ has a unique line $D$ of Whittaker functionals. So there is one and only one irreducible subquotient $\Pi$  with non-zero Whittaker functionals. As $\Pi$ is generic, $\Pi^{\s}$ is generic and by multiplicity one and the fact that $\ind_L^{G_E} \gamma$ is $\s$-stable we get $\Pi$ is $\s$-stable.

Set $\Pi=U/V$ with $(U,V)$ maximal as in the Appendix. Then 
$I_{\gamma,\s}(\ind_L^{G_E}\gamma)=\s^{-1} \circ \ind_L^{G_E} I_\gamma$ stabilizes $U$ and $V$, and induces by quotient $I_{\gamma,\s}(\Pi)$ on $\Pi$ (see the end of the Appendix).
 
If $\Lambda$ is a non-zero Whittaker functional in $D$, then it induces by restriction a non-zero Whittaker functional $\Lambda_U$ on $U$. 
By the same proof as in the previous Lemma \ref{gengen}, $I_{\gamma,\s}(\ind_L^{G_E}\gamma)$ fixes $\Lambda$, so its restriction to $U$ fixes $\Lambda_U$. This shows that
$I_{\gamma,\s}(\Pi)$ satisfies the definition of the \nio\ of $\Pi$. \qed

\ \\
{\bf Example.} Take for $\Pi$ the trivial character of $G_E$. When $n=1$, the trivial character is generic and $I_{\Pi}^{gen}$ is identity. 
For general $n$, $\Pi$ is isomorphic to $Lg(\Pi_1,\Pi_2,...,\Pi_n)$ where
$\Pi_i=\nu^{\frac{n-i}{2}}$, and each operator $I_{\Pi_i}^{gen}$ is the identity. 
Then the $\s$-intertwining operator on the space $V$ of $\Pi_1\t\Pi_2\t ...\t\Pi_n$ obtained by the parabolic induction procedure is composition of functions with $\s^{-1}$.
The space of left invariant functionals on $V$ has dimension one: one can construct one by using a Haar measure on $G_E$ and a left Haar measure on its upper triangular subgroup; such a functional is obviously invariant under $\s$, so the \nio\ of $\Pi$ is the identity.\\

Let $\chi$ be any character of $E^\t$, invariant under $\s$, i.e. factorizing through the norm from $E$ to $F$. Then if $\Pi$ is isomorphic to $\Pi^\s$, $\chi\Pi$ is also isomorphic to $(\chi\Pi)^\s$ and the \nio s are the same ($I_{\chi\Pi}=I_{\Pi}$, they act on the same space). From the example above, if $\Pi$ is any character of $G_E$ invariant under $\s$, then the \nio\ $I_\Pi$ is the identity.\\
\ \\
{\bf Remark.} Assume $(\Pi,V)$ is irreducible, $\s$-stable and unitary, and $(\ ,\ )$ is a $\Pi$-stable scalar product on $V$. A consequence of 
Schur's lemma is that every $\Pi$-stable scalar product on $V$ is proportional to $(\ ,\ )$. Obviously, the proportionality constant is real positive.
Set now $<v,v'>:=(I^{-1}_\Pi(v),I^{-1}_\Pi(v'))$. Then $<\ ,\ >$ is a scalar product, and is $\Pi^\s$-stable. But then it is also $\Pi$-stable. 
We have $<\ ,\, >=a(\ ,\ )$, with  $a\in \r_+^*$. Because $I_\Pi^l$ is identity by canonicity of $\Pi$, $a^l=1$ so $a=1$. So $I_\Pi$ {\it is a unitary operator of} $(V,(\ ,\ ))$.\\

\subsection{The norm map}

If $x\in\ge$, we write $Nx$
for the element $xx^\s...x^{\s^{l-1}}$ of $\ge$; it is called the {\bf
norm} of $x$.

Two elements $g,h\in\ge$ are called $\s${\bf -conjugate} if
$g=x^{-1}hx^\s$ for some $x\in\ge$.

\begin{lemme} {\rm (Lemma 1.1, page 3 [AC])}
{\rm (i)} If $x\in \ge$, $Nx$ is conjugate in $\ge$ to an element $y$ of
$\gf$;  $y$ is uniquely defined modulo conjugation in $\gf$.

{\rm (ii)} If $Nx$ and $Ny$ are conjugate in $\ge$, then $x$ and $y$ are
$\s$-conjugate.
\end{lemme}

So the norm map induces an injection from the set of $\s$-conjugacy classes
in $\ge$ into the set of conjugacy classes in $\gf$. If $x\in G_E$ (or if $x$ is a $\s$-conjugacy class in $G_E$) we will write
$\nn x$ for the associated conjugacy class in $\gf$.

If $P$ is a polynomial of degree $d$ with coefficients in $F$ we say $P$ is {\bf separable} if $P$ has $d$ distinct roots in an algebraic closure of $F$. If $y\in G_F$, we say $y$ is {\bf regular semisimple} if its characteristic polynomial is separable, and {\bf elliptic} if its characteristic polynomial is irreducible. A conjugacy class $y$ in $G_F$ is called regular semisimple if it contains a regular semisimple element (then all the elements in $y$ are regular semisimple, and regular semisimple classes in $G_F$ are parametrized by separable monic polynomials of degree $n$ with coefficients in $F$).
 
We say $x\in \ge$ is $\s${\bf -regular semisimple} if the class $\nn x$ is regular semisimple.

\subsection{The Shintani relation}\label{shintanilift}

Let $H(G_F)$ (resp. $H(G_E)$) be the Hecke algebra of complex functions locally constant with compact support on $G_F$ (resp. $G_E$). 
If $\pi$ is a finite length representation of $\gf$, the character of $\pi$ is the distribution on $\gf$ given by
$f\mapsto \tr(\pi(f))$ for $f\in H(\gf)$ (as always, $\pi(f):=\int_{G_F}f(g)\pi(g)dg$).
By [H-C], Theorem 1, or [DS] Theorem 16.3, there exists a locally integrable function $\chi_{\pi}$ on $G_F$ locally constant on the open set of regular semisimple elements, such that 
$\tr(\pi(f))=\int_{G_F} \chi_{\pi}(g) f(g) dg$ for $f\in H(\gf)$. 
Moreover, $\chi_\pi$ is constant on regular semisimple conjugacy classes.
Let $\Pi$ be an irreducible $\s$-stable representation of $\ge$. 
The twisted character of $\Pi$ is by definition the distribution on $\ge$ given by
$f\mapsto \tr(\Pi(f)\circ I_\Pi)$ for $f\in H(\ge)$. 
The operator $\Pi(f)$ is of finite rank, so $\tr(\Pi(f)\circ I_\Pi)=\tr(I_\Pi\circ \Pi(f))$ and we will use here one formula or the other\footnote{In [AC], the twisted character is defined, like in [La], with $I_\Pi$ on the right. However, the trace formula is written in [La] with the operator $M$ on the right and in [AC] it is written with the operator $M$ on the left, which makes $I_\Pi$ appear on the left in the trace formula of [AC] which we use here. But this has no influence on the local-global comparison because of the local equality of traces.}.
 
By [AC] Proposition 2.2, the twisted character is given by a locally integrable function $\chi_{\Pi,\s}$, locally constant on the open set of $\s$-regular
semisimple elements, i.e. $\tr(\Pi(f)\circ I_\Pi)=\int_{G_E} \chi_{\Pi,\s}(g)f(g)dg$. The function $\chi_{\Pi,\s}$ is constant on $\s$-regular conjugacy classes.

Let $\pi$ be an irreducible admissible representation
of $\gf$ and $\Pi$ an irreducible admissible representation of $\ge$. 
We say $\Pi$ and $\pi$ {\bf verify the Shintani relation}  if $\Pi$ is $\s$-stable and for all
$g\in \ge$ regular semisimple we have

$$\chi_{\Pi,\s}(g)=\chi_\pi(\nn g).$$

Note that this depends only on the isomorphism classes of $\pi$ and $\Pi$.

\subsection{Matching functions}\label{assocfunctions}

Let $T,T'$ be two maximal tori of $G_F$. Assume $T$ and $T'$ are conjugate, $T'=xTx^{-1}$ for some $x\in G_F$. A Haar measure $dt$ on $T$ induces then by conjugation a Haar measure $dt'$ on $T'$, which does not, actually, depend on the choice of $x$, as every continuous automorphism of $T$ is measure preserving. We say the measures $dt'$ and $dt$ are conjugate. On maximal tori of $G_F$ we fix Haar measures such that when two maximal tori are conjugate, the measures are conjugate. 

For $g\in G_F$, we let $G_g$ be the centralizer of $g$ in $G_F$. When $g$ is regular semisimple, $G_g$ is a maximal torus.

For  $g\in G_E$, let $G_{g,\s}$ be the twisted centralizer of $g$, namely the set of elements $y\in G_E$ such that $y^{-1}gy^\s=g$. 
Let $g\in G_E$ be $\s$-regular semisimple. 
Then, for every $\gamma\in \nn(g)$, there is a canonical $F$-isomorphism from $G_{g,\s}$ onto $G_{\gamma}$ ([AC], page 20). 
We fix on $G_{g,\s}$ the Haar measure defined by the pullback of the Haar measure chosen on $G_\gamma$. We define orbital integrals $\Phi$ and $\Phi_\s$ with respect to these measures, as in [AC], page 20: for $f\in H(G_F)$, $\gamma\in G_F$ regular semisimple,

$$\Phi(f,\gamma):=\int_{G_\gamma\bc G_F} f(x^{-1}\gamma x) dx$$
for the quotient measure and, for $\phi\in H(G_E)$, $g\in G_E$ $\s$-regular semisimple

$$\Phi_\s(\phi,g):=\int_{G_{\s,g}\bc G_E} \phi(x^{-1}gx^\s) dx$$
for the quotient measure.\\
\ \\
Then:

\begin{prop}{\rm (Proposition 3.1, page 20 [AC])}\label{intorb}

(a) If $\phi\in H(\ge)$ there exists $f\in H(\gf)$ such that, for all
regular semisimple $\gamma\in\gf$,

\ \ \ \ (i) $ \Phi(f,\gamma)=0$ if $\gamma$ is not conjugate to a norm and

\ \ \ \ (ii) $ \Phi(f,\gamma)=\Phi_\s (\phi,g)$ if $\gamma\in\nn g$.

(b) Given $f\in H(\gf)$ such that $ \Phi(f,\gamma)=0$ if $\gamma$ is not conjugate to a
norm, there exists $\phi\in H(\ge)$ such that

$$\Phi_\s(\phi,g)=\Phi(f,\nn g)$$
for all $g\in\ge$.
\end{prop}

If  $\phi$ and $f$ satisfy conditions (i) and (ii), we say $\phi$ and $f$ {\bf match}. We sometimes simply write $\phi \lra f$ to indicate that $\phi$ and $f$ match. The Weyl integration formula shows that $\Pi$ and $\pi$ verify the Shintani relation if and only if $\tr I_\Pi \Pi(\phi)=\tr\pi(f)$ whenever $f\lra \phi$. Then $\pi$ determines $\Pi$ up to isomorphism.\\
\ \\
{\bf Remark.} In [AC] the terminology ``$\phi$ and $f$ are associated" is used instead of ``$\phi$ and $f$ match". Here we made the choice to use ``match" in the local setting and ``associate" in the global setting where the nature of the definition is somehow different. (The same for ``square integrable" versus ``discrete series".)\\

Let $H^0(G_F)$ (resp. $H^0(G_E)$) be the sub-algebra of $H(G_F)$ (resp. $H(G_E)$) consisting of {\bf spherical} functions, i.e. functions which are left and right invariant by $K_F$ (resp. $K_E$). When $E/F$ is unramified, Arthur and Clozel define, using Satake parameters, an algebra morphism $b:H^0(G_E)\to H^0(G_F)$. They show {\rm ([AC] Theorem 4.5, page 39)} that $b(\phi)\lra  \phi$ for all $\phi\in H^0(G_E)$. This is the {\it fundamental lemma}. We will use the map $b$ later without explaining how it is constructed.

\subsection{Archimedean case} \label{archim}

Let us now consider briefly the Archimedean case, where $F\simeq \r$ or $\cc$.

The Langlands classification for $G_F=GL_n(F)$ has exactly the same statement as in the $p$-adic case above ([BW] Section IV, Theorem 4.11). Also, an irreducible unitary representation of a Levi subgroup of $G_F$ gives by normalized parabolic induction an irreducible unitary representation of $G_F$ ([Bar], and implicit in [Vo2]). The classification of irreducible unitary representations of $G_F$ is also the same as in the $p$-adic case ([Ta2]).

There are differences however between the Archimedean and the $p$-adic case. Firstly square integrable representations exist only when $n=1$ for $F\simeq \cc$, and $n=1$ or $n=2$ for $F\simeq \r$. Secondly there is some subtlety for generic representations. We explain it for an irreducible unitary representation $\Pi$ on a Hilbert space $V$, which is our case of interest. The subspace $V^{sm}$ of smooth vectors in $V$ carries a natural Fr\'echet space topology. On $V^{sm}$ we consider {\it continuous} linear functionals $\l$ such that 
$$\l(\pi(u)v)=\Theta_\psi(u)v$$
for all $u\in U_F$ and all $v\in V^{sm}$. The space of such functionals has dimension $0$ or $1$ ([Sh]) and we say that $\Pi$ is generic if it is $1$. Note, however, that on the $(\mathcal{G},K)$-module $V^\infty$ attached to $V$, there might be more functionals of this type ([Ko]). In any case, the classification of generic irreducible representations for $G_F$ is the same as in the $p$-adic case ([Vo], cf. the explanations in [He2], Section 2). Let now $E/F$ be the Archimedean extension $\cc/\r$ (up to isomorphism). Let $(\Pi,V)$ be an irreducible unitary representation of  $G_E$, invariant under $\s$ (i.e. under the conjugation of $E$ over $F$). As in the $p$-adic case, we construct, when $\Pi$ is generic, a \nio\ $I_\Pi^{gen}$ on $V^{sm}$. 

As in the Remark in Section \ref{niosection}, $I_\Pi^{gen}$ is unitary; it extends to a unitary operator on $V$ and
stabilizes $V^\infty$; we again write $I_\Pi^{gen}$ for the extended operator on $V$ and for the intertwining operator on the $(\mathcal{G},K)$-module $V^\infty$ induced by restriction. 
That is valid, in particular, when $\Pi$ is tempered. We then use the Langlands classification exactly as in the $p$-adic case to construct a  \nio\ $I_\Pi$ of $\Pi$ when $\Pi$ is not generic. When $\Pi$ itself is generic, then $I_\Pi=I_\Pi^{gen}$. To prove this, instead of relying on the results of Jacquet and Shalika used in the $p$-adic case, we use [Wa] 15.6.7, cf. the comments on the page 110 in [JS1]. In the same manner, we get a statement analogous to Proposition \ref{transitivity} in the Archimedean case, using [Wa] 15.6.7.

The norm map is defined exactly as in the $p$-adic case. Characters and twisted characters of representations also exist and the Shintani relation is defined via the same character identity as in section \ref{shintanilift}. Matching functions are defined in the same way and [AC] 7.3 gives the proposition analogous to \ref{intorb} for smooth $K$-finite functions with compact support.
Similarly, $\Pi$ and $\pi$ verify the Shintani relation
if and only if $\tr I_\Pi \circ \Pi(f)=\tr\pi(\phi)$ whenever we take such matching functions $f\lra \phi$. 
We also have $\tr  \Pi(f)\circ I_\Pi =\tr\pi(\phi)$, because for smooth $K$-finite functions with compact support $\Pi(f)$ is of finite rank.
We define Shintani lift (or simply lift) exactly as in \ref{shintanilift}.

\section{Results (local)} \label{resultslocal}

In this chapter \ref{resultslocal}, except in the last few lines, $E/F$ is a cyclic extension of $p$-adic fields of degree $l$ and that we have fixed a generator $\s$ of $Gal(E/F)$. Let $X(E/F)$ be the group of characters of $F^\t$ trivial on $N_{E/F}(E^\t)$. By class field theory, $X(E/F)$ is dual to $Gal(E/F)$, hence is cyclic of order $l$. We fix a generator $\chi$ of this group $X(E/F)$. If $\pi$ is a representation of $G_F$, we set $X_\pi$ for the set of isomorphism classes $\phi\pi$ where $\phi$ runs over $X(E/F)$. We let $m(\pi)$ be the cardinality of $X_\pi$ (it is the smallest positive integer $m$ such that $\pi\simeq \chi^m\pi$). It is clear that $m(\pi)$ divides $l$. If $\Pi$ is a representation of $G_E$, we let $X_\Pi$ be the set of isomorphism classes ${\Pi^x}$ where $x$ runs over $Gal(E/F)$ and $r(\Pi)$ be the cardinality of $X_\Pi$. Then $r(\Pi)$ is the smallest positive integer $r$ such that $\Pi^{\s^r}\simeq \Pi$. It is obvious that $r(\Pi)$ divides $l$.\\
\ \\
{\bf Theorem A} is proved in [AC], Theorem 6.2, Proposition 6.6, Lemma 6.10.\\
\ \\
{\bf Theorem A.}

{\it {\rm{(a)} \rm{(i)}} Let $\rho$ be a cuspidal representation of $GL_n(F)$.  Set $r:=l/m(\rho)$. Then there exists a $\s^{r}$-stable cuspidal representation $\rho'$ of $GL_{\frac{n}{r}}(E)$ such that 
$$\rho_{E}:=\rho'\tt  \rho'^\s \tt \rho'^{\s^2}\tt ...\tt  \rho'^{\s^{r-1}}$$ 
and  $\rho$ verify the Shintani relation. 
The representations $\pi$ of $GL_n(F)$ such that $\rho_{E}$ and $\pi$ satisfy the Shintani relation are up to isomorphism the elements of $X_{\rho}$. 
One has $r(\rho')=r$.

{\rm{(ii)}} Let 
$$\Pi=\pi'\tt  {\pi'}^\s \tt...\tt  {\pi'}^{\s^{r-1}}$$ 
be a representation of  $GL_{n}(E)$, where $\pi'$ is a cuspidal representation of  $GL_{\frac{n}{r}}(E)$ such that $r(\pi')=r$. Then there exists a cuspidal representation $\rho$ of $GL_n(F)$ such that $\Pi$ and $\rho$ verify the Shintani relation.  Moreover, one has $m(\rho)=l/r$. 

{{\rm (b)}} The same statement is true if we replace ``cuspidal" with ``square integrable" in {\rm{(a)}}. Moreover we have

{\rm{(iii)}} if $\d=Z(\rho,k)$ is a square integrable representation of $GL_n(F)$, then $m(\rho)=m(\d)$. If $\rho$ and $\rho'$ are linked by the relation of {\rm (a) (i)}, set $\d'=Z(\rho',k)$. Then 
$$\d_{E}:=\d'\tt  \d'^\s \tt \d'^{\s^2}\tt ...\tt  \d'^{\s^{r-1}}$$ 
and $\d$ verify the Shintani relation.

{{\rm (c) (i)}} Let $\tau$ be a tempered representation of $GL_n(F)$. There is a (unique) tempered $\s$-stable representation $\tau_{E}$ of $GL_n(E)$ such that $\tau_E$ and $\tau$ verify the Shintani relation.

{{\rm (ii)}} If $T$ is an irreducible $\s$-stable tempered representation of $GL_n(E)$, then there is a tempered representation $\tau$ of $GL_n(F)$ such that $T$ and $\tau$ verify the Shintani relation.}\\
\ \\

The proof of this Theorem in [AC] is written only for prime $l$ (so that $m(\rho)$ is always $1$ or $l$) but works in general. For a tempered representation $\tau$ of $GL_n(F)$, we say that the representation $\tau_E$ provided by the {\bf Theorem A} (c) (i) is {\bf the base change}, and also {\bf the Shintani lift}, of $\tau$. Using the Langlands classification, Arthur and Clozel define an abstract {\bf base change} for any irreducible representation $\pi$ of $G_F$ to a $\s$-stable representation $\pi_{E}$ of $G_E$. If $\pi=Lg(\tau_1,\tau_2,...,\tau_k)$ for essentially tempered representations $\tau_1,\tau_2,...,\tau_k$, then $\pi_{E}$ is defined by 
$\pi_{E}=Lg(\tau_{E,1}, \tau_{E,2},...,\tau_{E,k})$, where $\tau_{E,i}$ is the Shintani lift (from (c) (i)) of $\tau_i$ for $1\leq i\leq k$.
Here we say that $\pi_E$ (the base change of $\pi$) is {\bf a Shintani lift} of $\pi$ if $\pi_E$ and $\pi$ satisfy the Shintani relation. Thus it is automatic for tempered $\pi$, but that is NOT always true in general, as explained in our example after Proposition \ref{liftprod}. When $\pi_E$ is a Shintani lift of $\pi$, we also say that $\pi$ has a Shintani lift. Our goal is to provide many classes of representations of $GL_n(F)$ which do have a Shintani lift.\\
\ \\
{\bf Theorem B.}
{\it Every irreducible spherical unitary representation of $GL_n(F)$ has a Shintani lift which is an irreducible spherical unitary representation of $GL_n(E)$.
Every irreducible spherical unitary representation of $GL_n(E)$ is the Shintani lift of an irreducible spherical unitary representation of $GL_n(F)$.}
\ \\



That is proved in [AC] (I. 4, III.1) using base change for spherical functions. We give here a straightforward proof (see section \ref{proofs} for the proof and an explicit form of this Shintani lift).\\
\ \\
The original local results of this paper are the following. If $\pi$ is a smooth representation of $G_F$ (or $G_E$), we say $\pi$ is {\bf elliptic} if $\pi$ is irreducible and the character $\chi_\pi$ of $\pi$ does not vanish identically on the elliptic set, i.e. $\tr\pi(f)$ is not identically null on functions $f$ with support in the elliptic set.\\
\ \\
{\bf Theorem C.}  

{\it {\rm{(a)}} Let $\xi$ be an elliptic representation of $GL_n(F)$.  Set $r:=l/m(\xi)$. Then there exists a $\s^{r}$-stable elliptic representation $\xi'$ of $GL_{\frac{n}{r}}(E)$ such that 
$$\xi_{E}:=\xi'\tt  \xi'^\s \tt \xi'^{\s^2}\tt ...\tt  \xi'^{\s^{r-1}}$$ 
is a Shintani lift of $\xi$. All the representations of $GL_n(F)$ with Shintani lift $\xi_{E}$  are up to isomorphism the elements of $X_{\xi}$. One has $r(\xi_{E})=r$.

{\rm{(b)}} Let 
$$\Xi=\xi'\tt  {\xi'}^\s \tt...\tt  {\xi'}^{\s^{r-1}}$$ 
be a representation of  $GL_{n}(E)$, where $\xi'$ is an elliptic representation of  $GL_{\frac{n}{r}}(E)$ such that $r(\xi')=r$. Then $\Xi$ is the Shintani lift of some elliptic representation $\xi$ of $GL_n(F)$. Moreover, one has $m(\xi)=l/r$ and the representations with Shintani lift $\Xi$ are up to isomorphism the elements of $X_{\xi}$.}


\ \\
{\bf Theorem D.}

{\it {\rm{(a)}} Let $\d$ be a square integrable representation of $GL_n(F)$ and let
$$\d_{E}=\d'\tt  \d'^{\s} \tt \d'^{\s^2}\tt ...\tt  \d'^{\s^{r-1}}$$ 
be the Shintani lift of $\d$ to $GL_n(E)$ like in {\bf Theorem A}.  Let $k$ be a positive integer and set $\pi:=u(\d,k)$. Then $m(\pi)=m(\d)=l/r$ and 
$$\pi_{E}:=u\t u^\s\t u^{\s^2}\t ...\t u^{\s^{r-1}},$$
where $u:=u(\d',k)$, is a Shintani lift of $\pi$. 

{\rm{(b)}} Every unitary irreducible representation $\tau$ of $GL_n(F)$ has a Shintani lift $T$ to $GL_n(E)$. Every unitary irreducible $\s$-stable representation $T$ of $GL_n(E)$ is the Shintani lift of some unitary irreducible representation $\tau$ of $GL_n(F)$. The representation $\tau$ is generic if and only if $T$ is generic.}\\
\ \\
\ \\
{\bf The case $F\simeq \r$.} In the case $F\simeq\r$ and $E\simeq \cc$, Theorems {\bf A}, {\bf C} and {\bf D} are still true (for {\bf B}, the notion of spherical representation does not make sense). In this case $l=2$ and so $m,r\in\{1,2\}$. Moreover, cuspidal representations exist only in the case $n=1$, square integrable representations exist only if $n\in\{1,2\}$, elliptic elements exist in the case $n\in\{1,2\}$ and then the elliptic representations are the essentially square integrable ones and the finite dimensional ones.

\section{Notation and basic facts (global)}\label{notglobal}

Let $\F$ be a number field and $n$ a positive integer. For every place $v$ of $\F$ we let $\F_v$  be 
the $v$-adic completion of $\F$. For each place $v$ we fix an absolute value $|\ \ |_v$ on $\F_v$, normalized as before if $v$ is finite, and equal to the usual absolute value (resp. the square of the modulus) if $\F_v\simeq\r$ (resp. $\F_v\simeq \cc$).
If $v$ is finite, let $O_{\F,v}$ be the ring of integers of $\F_v$. Let $\aa_{\F}$ be the ring of ad\`eles of $\F$. 
For every place $v$ of $\F$ we denote $G_{\F,v}$ the group $GL_n(\F_v)$. 
If $v$ is finite we let $K_{\F,v}$ be the maximal compact subgroup $GL_n(O_{\F,v})$ of $G_{\F,v}$. 
We then endow $G_{\F,v}$ with the Haar measure such that the volume of $K_{\F,v}$ is one. If $v$ is infinite, we fix a Haar measure $dg_v$ on $G_{\F,v}$.

Let
$GL_n(\aa_{\F})$ be the ad\`ele group of $GL_n(\F)$ with respect to the $K_{\F,v}$. 
We denote $|\ \ |_{\aa_\F}$ the absolute value on $\aa_\F$ which is the product of the local absolute values. 
If we see $\F$ as a subring of $\aa_\F$ (by diagonal embedding), then elements of $\F^\t$ have absolute value $1$. We endow $GL_n(\aa_{\F})$ with the Haar measure $dg$ which is the product of local Haar measures $dg_v$.

We consider $GL_n(\F)$ as a subgroup of $GL_n(\aa_{\F})$ via the diagonal embedding. Then $GL_n(\F)$ is a discrete subgroup of the locally compact group $GL_n(\aa_\F)$.

Let $Z$ be the center of $GL_n$, $Z(\F)$ the center of $GL_n(\F)$. For every place $v$, let $Z_{\F,v}$ be the center of $G_{\F,v}$. 
For every finite place $v$ of $\F$, let $dz_v$ be a Haar measure on $Z_{\F,v}$ such that the volume of $Z_{\F,v}\cap K_{\F,v}$ is one.   For infinite places $v$, fix a Haar measure $dz_v$ on $Z_{\F,v}$.
The center $Z(\aa_{\F})$ of $GL_n(\aa_{\F})$ is the restricted product of the $Z_{\F,v}$ with respect to the $Z_{\F,v}\cap K_{\F,v}$. 
On $Z(\aa_{\F})$ we fix the Haar measure $dz$ which is the product of the measures $dz_v$. 
On $Z(\aa_{\F})\bc GL_n(\aa_{\F})$ we consider the quotient measure $dz\bc dg$. 
As $GL_n(\F)\cap Z(\aa_{\F})\bc GL_n(\F)$ is a discrete subgroup of $Z(\aa_{\F})\bc
GL_n(\aa_{\F})$, on the quotient space $Z(\aa_{\F})GL_n(\F)\bc GL_n(\aa_{\F})$ we have a well-defined measure $d\bar{g}$ coming from $dz\bc dg$. 
The measure of the whole space $Z(\aa_{\F})GL_n(\F)\bc GL_n(\aa_{\F})$ is finite.

Let $\o$ be a unitary character of $Z(\aa_\F)$, trivial on $Z(\F)$. We write $\lg$ for 
the space of classes (modulo zero measure sets) of functions $f$ defined on
$GL_n(\aa_{\F})$ with values in $\cc$ such that:

i) $f$ is left invariant under $GL_n(\F)$,

ii) $f$ satisfies $f(zg)=\o(z)f(g)$ for all $z\in Z(\aa_{\F})$ and almost
all $g\in GL_n(\aa_{\F})$,

iii) $|f|^2$ is integrable over $Z(\aa_{\F})GL_n(\F)\bc GL_n(\aa_{\F})$.\\
The scalar product on $\lg$ is the standard one, given by
$$(f,h)=\int_{Z(\aa_{\F})GL_n(\F)\bc GL_n(\aa_{\F})}\overline{f(\bar{g})}h(\bar{g})d\bar{g},$$
where $\overline{f(\bar{g})}$ is the complex conjugate of $f(\bar{g})$.

We consider the representation $R_\o$ of $GL_n(\aa_{\F})$ by right translations on this space. We call a representation of 
$GL_n(\aa_{\F})$ a {\bf discrete series} if it is an irreducible subrepresentation (irreducible stable Hilbert subspace) of such a representation 
$R_\o$ for any unitary character $\o$ of $Z(\aa_{\F})$ trivial on $Z(\F)$. The underlying automorphic representation ([BJ]) is irreducible. We call such an automorphic representation {\bf automorphic discrete series}. They can be either cuspidal or residual. Recall that a discrete series $\Pi$ is a completed Hilbert tensor product of local Hilbert unitary representations $\Pi_v$, and the underlying automorphic representation $\Pi^{\infty}$ is the restricted tensor product of representations $\Pi_v^\infty$, where $\Pi_v^\infty$ is the $({\mathcal G},K)$-module of $Z({\mathcal G})$-finite and $K$-finite vectors in $\Pi_v$ if $v$ is infinite, and is the space of smooth vectors in $\Pi_v$ if $v$ is finite. 

Let ${\mathcal L}(\F)$ be the set of Levi subgroups of $GL_{n}(\F)$ given by block matrices as in the local setting. For $k$ dividing $n$, let $L_{\F,k}$ be the Levi 
subgroup of $GL_n(\F)$ given by block diagonal matrices with $k$ blocks of the same size (which is $n/k$). 
By abuse of notation, if $L\in {\mathcal L}(\F)$, we will denote by the same letter $L$ the corresponding Levi subgroup of $GL_n(\aa_{\F})$ when no confusion may occur.

Moeglin and Waldspurger gave ([MW1]) the classification of automorphic discrete series in terms of automorphic cuspidal representations.
If $k|n$ and $\rho$ is an automorphic cuspidal representation of $GL_{n/k}(\aa_\F)$, then the representation  
$\ind_{L_k}^{GL_n(\aa_\F)} (\nu^{\frac{k-1}{2}}\rho\otimes \nu^{\frac{k-3}{2}}\rho\otimes \nu^{\frac{k-5}{2}}\rho\otimes ... \otimes\nu^{-\frac{k-1}{2}}\rho)$
has a unique irreducible quotient $MW(\rho,k)$ which is an automorphic discrete series of $GL_n(\aa_\F)$ (see [La2] for the definition of parabolic induction and quotient for global representations). Given an automorphic discrete series $\d$ of $GL_n(\aa_\F)$ there exist $k|n$ and an automorphic cuspidal representation $\rho$ of $GL_{n/k}(\aa_\F)$ such that $\d$ is isomorphic to $MW(\rho,k)$. This is the Th\'eor\`eme on page 606, [MW1]. Let us show that $\d$ actually determines $k$ and the class of $\rho$.  
We show it together with the {\it strong multiplicity one theorem for discrete series}.\\
\ \\
{\bf Fact.} {\it If $\d=MW(\rho,k)$ and $\d'=MW(\rho',k')$ are two automorphic discrete series of $GL_n(\aa_\F)$ such that $\d_v \simeq \d'_v$ for almost all places $v$ of $\F$, then $\d_v \simeq \d'_v$ for all places $v$ of $\F$, i.e. $\d\simeq \d'$. Moreover, $k=k'$ and $\rho\simeq \rho'$.}\\

This result is known to hold for cuspidal representations, i.e. the case $k=k'=1$ ([Sh], [P-S]). Let $v$ be a finite place of $\F$.
If $\d=MW(\rho,k)$, then by construction (see below some detail in the proof of Proposition \ref{opintMW}).
$\d_v=u(\rho_v,k)$, where $\rho_v$ is the local component of $\rho$ at the place $v$, which is known to be generic by [Sh]. 
If $\d'=MW(\rho',k')$, then $\d'_v=u(\rho'_v,k')$, where $\rho'_v$ is the local component of $\rho'$ at the place $v$. The Tadi\'c classification of unitary representations and the fact that a unitary generic representation is an irreducible product of essentially square integrable representations implies that $u(\rho_v,k)\simeq u(\rho'_v,k')$ implies $k=k'$ and $\rho_v\simeq \rho'_v$. So the strong multiplicity one theorem for automorphic discrete series follows from the strong multiplicity one theorem applied to the cuspidal representations $\rho$ and $\rho'$.\\

Let $\d=MW(\rho,k)$ be a discrete series of $GL_n(\aa_\F)$. We then have the {\it multiplicity one theorem for automorphic discrete series}:\\
\ \\
{\it $\d$ appears with multiplicity one in $\lg$.}\\

This is the Th\'eor\`eme of [MW1] page 606, combined with the fact that $\d$ determines $k$ and $\rho$.\\
\ \\

Let $\E$ be a cyclic extension of $\F$. Adopt the same notation for $\E$ as for $\F$. 

Put $l=[\E :\F ]$ and choose a generator $\s$ of the Galois group $Gal(\E/\F)$ of $\E/\F$.

Let $v$ be a place of $\F$ and set $\E_v=\E\otimes_\F \F_v$. 
Then $\E_v$ decomposes naturally as the product $\prod_{w|v} \E_w$, where $w$ runs through the places of $\E$ above $v$, and $Gal(\E/\F)$ acts $\F_v$-linearly on $\E_v$, permuting the $\E_w$. More precisely, the extensions $\E_w/\F_v$ all have the same degree $d$ and $l=da$, where $a$ is the number of places $w$ above $v$; the stabilizer of $\E_w$ in $Gal(\E/\F)$ is generated by $\s^a$. 
We can choose notation so that $\E_v$ is the product $\E_1\t\E_2\t ...\t\E_a$ (of extensions of $\F_v$), with $\s(\E_i)=\E_{i+1}$ for $i=1,2,...,a-1$, $\s(\E_a)=\E_1$. 
An irreducible smooth representation $\Pi$ of $GL_n(\E_v)=GL_n(\E_1)\t GL_n(\E_2)\t ...\t GL_n(\E_a)$ then decomposes as a tensor product $\Pi_1\otimes \Pi_2\otimes ...\otimes \Pi_a$ and $\Pi$ is $\s$-stable if and only if $\Pi_1$ is $\s^a$-stable and $\Pi_i\circ\s^{i-1}\simeq \Pi_1$ for $i=1,2,...,a$.

If that is the case, we can take $\Pi_i=\Pi_1\circ\s^{1-i}$ and if $I_1$ is an isomorphism of $\Pi_1$ onto $\Pi_1^{\s^a}$, then $I:x_1\otimes x_2\otimes ...\otimes x_a\mapsto I_1x_a\otimes x_1\otimes...\otimes x_{a-1}$ is an isomorphism of $\Pi$ onto $\Pi^\s$. If $I_1$ is the normalized $\s^a$-intertwining operator of $\Pi_1$, we call $I$ the \nio\ of $\Pi$ (it is readily verified that it does not depend on the identifications $\Pi_i=\Pi_1\circ\s^{1-i}$).\\

The notions of  matching functions and Shintani lift readily extend to the case of cyclic $\F_v$-algebra $\E_v$ (see [AC] I, 5, and more explanation in [HL]). 
If $\Pi$ is an irreducible $\s$-stable representation of $GL_n(\E_v)$ written $\Pi=\Pi_1\otimes \Pi_2\otimes ...\otimes \Pi_a$ as above, then $\Pi$ is a Shintani lift of $\pi$ if and only if, for $i=1,2,...,a$, $\Pi_i$ is a Shintani lift of $\pi$, with respect to  
the generator $\s^a$ of $Gal(\E_i/ \F_v)$; as $\Pi$ is $\s$-stable, it is enough to verify that $\Pi_1$ is a Shintani lift of $\pi$ with respect to $\s^a$. 
Now $H(GL_n(\E_v))$ is identified with $H(GL_n(\E_1))\otimes H(GL_n(\E_2))\otimes ...\otimes H(GL_n(\E_a))$, 
and if $\phi_v=\phi_1\otimes \phi_2\otimes ...\otimes \phi_a\in H(GL_n(\E_v))$, 
we will view the $\phi_i$ as functions on $GL_n(\E_1)$ via the isomorphism $\s^{i-1}:\E_1\to \E_i$ and associate to $\phi_v$ the function 
$\tilde{\phi}_v:=\phi_1*\phi_2* ...*\phi_a\in H(GL_n(\E_1))$ (where $*$ is the convolution product). 
We then say $\phi_v$  and $f_v$ match if $\tilde{\phi}_v$ and $f_v$ match in the sense of section \ref{assocfunctions} (see [HL] I.2.5, I.2.9 for details). 
Notice that $\tr (I\circ \Pi(\phi_v))=\tr( (I_1\circ\Pi_1(\phi_1))\circ\Pi_1(\phi_2)\circ\Pi_1(\phi_3)...\circ\Pi_1(\phi_a))=\tr (I_1\circ\Pi_1(\tilde{\phi}_v))$. 
We extend the definition of spherical base change $b$ (section \ref{assocfunctions}) by setting $b(\phi_v)=b(\tilde{\phi}_v)$.\ \\

The Galois group $Gal(\E/\F)$ acts on $GL_n(\aa_\E)$, hence also on its unitary or admissible representations. As in the local setting, if $\Pi$ is such a representation, $\Pi^\s$ will be the representation on the same space given by $g\mapsto \Pi(\s(g))$. Also, $Gal(\E/\F)$ acts on functions from $GL_n(\aa_\E)$ to $\cc$. As in the local setting, if $\Pi$ is a representation  of $GL_n(\aa_\E)$ by right translations on some space $V$ of functions from $GL_n(\aa_\E)$ to $\cc$, we let $\s V=\{f\circ\s | f\in V\}$ an $\s\Pi$ be the representation of $GL_n(\aa_\E)$ by right translations in $\s V$. 
Note that the map $f\mapsto f\circ\s$, $V\to \s V$, is an isomorphism of $\Pi^\s$ onto $\s\Pi$.
This applies, for example, to the representation $R_\o$ on $\lge$ for any unitary character $\o$ of $\aa_\E^\t$ trivial on $\E^\t$. If $\o$ is invariant under $Gal(\E/\F)$, then $\s (\lge) = \lge$.

A discrete component of $\lge$ occurs with multiplicity one.
So, if $(\Pi,V)$ is such a component and $\Pi\simeq \Pi^\s$ then actually $V=\s V$ and $f\mapsto f\circ\s$ gives an isomorphism of $\Pi^\s$ onto $\s \Pi=\Pi$. 
The operator $\s^{-1}:f\mapsto f\circ \s^{-1}$ is called the ``physical" operator between $\Pi$ and $\Pi^\s$. But we also have a \nio\ $I_\Pi$ coming from the local setting, as we now explain. For each place $v$ of $\F$, let $\Pi_v$ be the local component of $\Pi$ at $v$. It is an irreducible unitary representation of $GL_n(\E_v)=\Pi_{w|v} GL_n(\E_w)$, well defined up to isomorphism. Choose an isomorphism $\iota$ of $\hat{\otimes}_v \Pi_v$ onto $\Pi$ where $\hat{\otimes}$ stands for the completed restricted tensor product (restricted with respect to the choice of a unitary spherical vector $e_v$ in $\Pi_v$ for almost all finite places $v$, see [Fl1]). As $\Pi$ is isomorphic to $\Pi^\s$, we also have $\Pi_v\simeq \Pi_v^\s$ for each place $v$ of $\F$, so we get an associated \nio\  $I_{\Pi_v}$; it is unitary and stabilizes $e_v$ for almost all $v$ (see above, Section \ref{niosection}). So ${\otimes}_v I_{\Pi_v}$ 
has a unique extension to a unitary operator $\hat{\otimes}_v I_{\Pi_v}$ on $\hat{\otimes}_v \Pi_v$ which transfers 
via $\iota$ to the operator $I_\Pi:\Pi\to \Pi^\s$, which does not depend on the choice of $\iota$ (nor on the choice of $e_v$'s).

\begin{prop}\label{opintMW}
Let $\Pi$ be a discrete series component of ${L^2(GL_n(\F)Z(\aa_\F)\bc GL_n(\aa_\F),\o)}$ for some unitary character $\o$ of $\aa_\E^\t$ trivial on $\E^\t$ and $\s$-invariant. Assume that $\Pi\simeq \Pi^\s$. Then $I_\Pi$ is equal to the physical operator $f\mapsto f\circ \s^{-1}$.
\end{prop}
\ \\
{\bf Proof.} The case where $\Pi$ is cuspidal is known: indeed it is present behind the scene in [AC]. We give an argument, based on [Sh].

Let $\Pi^{sm}$ be the image by $\iota$ of $\otimes_{v|\infty}\Pi_v^{sm} \otimes_{v\text{ finite}} \Pi_v^{\infty}$. 
Here, as above in section \ref{archim}, $\Pi_v^{sm}$ denotes the space of smooth vectors in $\Pi_v$ if $v$ is infinite, and $\Pi_v^\infty$ is the space of vectors in $\Pi_v$ with open stabilizer if $v$ is finite. 
On $\Pi^{sm}$ we have the global Whittaker functional $\Lambda: f\mapsto \int_{U(\E)\bc U(\aa_\E)} f(u) \Theta_{\psi_\E}(u) du$, where we have chosen a non-trivial character  $\psi$ of $\aa_{\F}/\F$ and put $\psi_\E=\psi\circ \tr_{\E/\F}$, defining $\Theta_{\psi_\F}$ as in the local case, and we have chosen a Haar measure $du$ on $U(\aa_\E)$. 
By (loc. cit.) $\Lambda$ is factorisable i.e. $\Lambda\circ \iota$ is the tensor product $\otimes \Lambda_v$ of local Whittaker functionals $\Lambda_v$ with respect to $\psi_{\E_v}=\psi_{\E_{|\E_v}}$. 
Since Haar measures on $U(\aa_\E)$ are $\s$-invariant and $\Theta_{\psi_\E}$ is  $\s$-invariant too, we have $\Lambda(f)=\Lambda(f\circ\s^{-1})$for $f$ in the space of $\Pi^{sm}$. 
On the other hand, the \nio\ $I_v$ on $\Pi_v^\infty$ or $\Pi_v^{sm}$ is the unique one such that $\Lambda_v\circ I_v=\Lambda_v$. So we get that $\otimes I_v$ transports to the physical operator on $\Pi^{sm}$ and, by density, on $\Pi$.

To treat the case of residual $\Pi$, we need to recall the construction by Jacquet ([Ja]) of such representations.\\

Jacquet's construction starts with a strict divisor $r$ of $n$, $n=ra$, and a cuspidal automorphic representation $\Sigma^\infty$ of $GL_r(\aa_\E)$ on a space, say $V^\infty$, of functions on $GL_r(\E)\bc GL_r(\aa_\E)$. 
Jacquet assumes in fact, in the ambiguous ``we may arrange" ([JS], p. 187, line -4) that the central character of $\Sigma^\infty$ is trivial on the subgroup $\r_+^\t$ of $\aa_\E^\t$, where $\aa_\E^\t$ is seen as the center of $GL_r(\aa_\E)$ and $\r_+^\t$ is embedded diagonally at the infinite places of $\E$. 
We assume that condition for the moment -- we shall say at the end what to do in the general case.

Write $M$ for the block diagonal Levi subgroup of $GL_{n}$ with blocks of size $r$ (so $M$ is the Levi subgroup $L_F$ in the above notation), $P$ for the corresponding upper triangular parabolic subgroup, and $U_P$ for the unipotent radical of $P$. Jacquet considers the automorphic representation of $GL_n(\aa_\E)$ automorphically induced from 
$\Sigma^\infty\otimes\Sigma^\infty\otimes...\otimes\Sigma^\infty$, 
which he realizes in a space $\mathcal{F}$ of functions $f$ on $GL_n(\aa_\E)$ which are $K$-finite (where $K$ is the product over all places $w$ of $\E$ of the usual maximal compact subgroup $K_w$ of $GL_n(\E_w)$) left-invariant under $U_P(\aa_\E)$ and $P(\E)$ and such that for each $k$ in $K$ the function $m\mapsto f(mk)$ on $M(\aa_\E)$ belongs to the space 
$V^\infty\otimes V^\infty\otimes ...\otimes V^\infty$ 
of automorphic forms on $M(\aa_\E)$. For $s\in \cc^a$ one forms an Eisenstein series 
$$E(g,s,f)=\sum_{\gamma\in P(\E)\bc G(\E)} f(\gamma g) \exp (<H_P(\gamma g),s+\rho>)$$

That Eisenstein series converges absolutely in an open cone in $\cc^a$ and extends to a meromorphic function of $s$. Jacquet specifies a meromorphic function $u(s)$ (independent of $f$ and whose exact value is irrelevant for us) such that $u(s) E(g,s,f)$ is holomorphic at $e:=(\frac{a-1}{2},\frac{a-1}{2}-1,...,\frac{1-a}{2})$ 
(entries decreasing by $1$) and such that, varying $f$, we get upon evaluation at $s=e$ the desired space of automorphic functions on $GL_n(\aa_\E)$, spanning the automorphic discrete series $\Pi^\infty$ attached to $\Sigma^\infty$. This process gives a $GL_n(\aa_\E)$-intertwining operator between the space ${\mathcal{F}}_e$ of functions $g\mapsto f(g)\exp(<H_P(g),e+\rho>)$ and $\Pi^{\infty}$. We write $Res$ for the resulting map $\mathcal{F}_e\to \Pi^\infty$.

For a place $w$ of $\E$, let $\Sigma_w$ be the local component of $\Sigma^\infty$ at $w$, $V_w$ the space of $\Sigma_w$, and form the representation 
$R_w:= \nu_{\E_w}^{\frac{a-1}{2}}\Sigma_w\t \nu_{\E_w}^{\frac{a-1}{2}-1}\Sigma_w\t ...\t \nu_{\E_w}^{\frac{1-a}{2}}\Sigma_w$ of $GL_n(\E_w)$. 
Jacquet identifies the local component $\Pi_w$ of $\Pi^\infty$ at $w$ as the Langlands quotient of $R_w$ (note that that quotient is the representation $u(\Sigma_w,a)$; we used that in the proof of the fact above). 
As [Ja] is very elliptic, we give enough detail to follow the action of the $\s$-operator. 

Let $J$ be an isomorphism of $\otimes_w \Sigma_w$ onto $\Sigma^\infty$, so that $J$ is a linear isomorphism of $\otimes_w V_w$ onto the space $V^\infty$ of $\Sigma^\infty$. An element $\phi_w$ of $R_w$ is a function $GL_n(\E_w)\to V_w\otimes V_w\otimes ...\otimes V_w$ satisfying some extra conditions. To an element $\phi=\otimes\phi_w$ of $\otimes_w R_w$ we associate a function 
$$\mathcal{F}_\phi:GL_n(\aa_\E)\to \cc$$ 
in the following way. First define the function $\Phi_\phi: GL_n(\aa_\E)\to V^\infty\otimes V^\infty\otimes...\otimes V^\infty$, 
$\Phi_\phi:g=(g_w)_w\mapsto (J\otimes J\otimes ...\otimes J) (\otimes_w\phi_w(g_w))$. Recalling that $V^\infty$ is actually made out of functions $GL_r(\aa_\E)\to \cc$, we can evaluate at $1\in GL_r(\aa_\E)$ a function in $V^\infty$, resulting in a linear map $V^\infty\otimes V^\infty\otimes ...\otimes V^\infty\to\cc$. Composing with $\Phi_\phi(g)$ yields the desired function $g\mapsto \mathcal{F}_\phi(g)$ (so, if $\phi=\phi_1\otimes\phi_2\otimes...\otimes\phi_a$, $\mathcal{F}_\phi(g)=\prod_{i=1}^aJ(\phi_i(g))(1)$ where the product is taken in $\cc$).

This process results in an isomorphism $\Phi$ of $\otimes_w R_w$ onto $\mathcal{F}_e$, as automorphic representation of $GL_n(\aa_\E)$. Composing with $Res$ gives a surjective $GL_n(\aa_\E)$-equivariant map $\otimes_w R_w\to \Pi^\infty$. The local component $\Pi_w$ of $\Pi^\infty$ at $w$ can then only be the Langlands quotient $\bar{R}_w$ of $R_w$, and the $GL_n(\aa_\E)$-equivariant map $\otimes_w R_w\to \Pi^\infty$ factors through an isomorphism $\otimes_w\bar{R}_w\simeq \Pi^\infty$.

\def\ii{{\mathcal I}}
Now assume that $\Pi$ is $\s$-stable. That happens exactly when $\Sigma$ is $\s$-stable, by the uniqueness result of Moeglin and Waldspurger. We clearly have 
$Res(f\circ\s^{-1})=Res(f)\circ\s^{-1}$ for $f\in \mathcal{F}_e$. For each place $w$ of $\E$, let $I_w$ be the \nio\  on $\Sigma_w$, and $\ii_w$ the $\s$-operator on $R_w$ obtained by the parabolic induction procedure from $I_w$ (using our previous notation, $\ii_w=I_\s(\Sigma_w)$). Then $\ii_w$ induces the \nio\  on $\bar{R}_w$. Thus, it is enough to prove that for $\phi=\otimes_w\phi_w$ in $\otimes_w R_w$, we have ${\mathcal F}_{\otimes_w \ii_w\phi_w}(g)={\mathcal F}_\phi(\s^{-1}(g))$ for $g\in GL_n(\aa_\E)$.

But $\ii_w\phi_w$ is the function $g_w\mapsto (I_w\otimes I_w\otimes...\otimes I_w)(\phi_w(\s^{-1}(g_w)))$, so we get 
$\Phi_{\otimes_w\ii_w\phi_w}(g)=(J\otimes J\otimes ...\otimes J)((\otimes_wI_w)\otimes (\otimes_wI_w) \otimes ...\otimes (\otimes_wI_w))
(\otimes_w \phi_w(\s^{-1}(g_w)))$. By the cuspidal case already treated, $J\circ(\otimes_w I_w)$ is $I\circ J$, where $I$ is the physical operator on $V^\infty$, so we get
$\Phi_{\otimes_w\ii_w\phi_w}(g)=(I\otimes I\otimes ...\otimes I)(\Phi_{\otimes\phi_w}(\s^{-1}(g))$.
Composing with evaluation in $1$, we find ${\mathcal F}_{\otimes_w\ii_w\phi_w}(g)= {\mathcal F}_{\otimes \phi_w}(\s^{-1} g)$ (since $(I\circ u)(1)=u(1)$ for $u\in V^\infty$), which is what we wanted. 

The result for $\Pi$ follows by density from the result on $\Pi^\infty$.

Finally, when the central character of $\Sigma$ is not trivial on $\r_+^\t$, there is some complex number $t$ such that the central character of $\Sigma'=|\ |_\E^{-t}\Sigma$ (realized in the space of functions $g\mapsto |\det g|_\E^{-t}u(g)$ for $u$ in $V^\infty$) is trivial on $\r_+^\t$. If $\Pi'$ is the discrete series attached to $\Sigma'$ by the above procedure, then $\Pi=|\ |_\E^t \Pi'$ is the one attached to $\Sigma$ -- it is realized in the space of functions $g\mapsto |\det|_\E^t\phi(g)$, $\phi$ in the space of $\Pi'$. 
From the case of $\Pi'$ established above, one deduces right away that the \nio\  on $\Pi$ is equal to the physical operator.\qed

\section{Results (global)}

Let $\E/\F$ be a cyclic extension of number fields of degree $l$ as before. Let $X(\E/\F)$ be the group of characters of  $\aa_{\F}^\tt$ trivial on $\F^\tt N_{\E/\F}(\aa_{\E}^\tt)$.
By class field theory, $X(\E/\F)$ is isomorphic to the dual of $Gal(\E/\F)$, hence is cyclic of order $l$. We fix a generator $\chi$ of this group. If $\pi$ is a discrete series of $G(\aa_\F)$, we set $X_\pi$ for the set of isomorphism classes $\phi\pi$ where $\phi$ runs over $X(\E/\F)$. We let $m(\pi)$ be the cardinality of $X_\pi$. If $\Pi$ is a discrete series of $GL_n(\aa_\E)$, we let $X_\Pi$ be the set of isomorphism classes ${\Pi^x}$ where $x$ runs over $Gal(\E/\F)$ and we let $r(\Pi)$ be the cardinality of $X_\Pi$. Then $r(\Pi)$ is the smallest positive integer $r$ such that $\Pi^{\s^r}\simeq \Pi$. 

If $\pi$ is an automorphic discrete series of $GL_n(\aa_\F)$ and $\Pi$ is a $\s$-stable irreducible automorphic representation of $GL_n(\aa_\E)$, we say $\Pi$ is a Shintani lift of $\pi$ if, for every place $v$ of $\F$, $\Pi_v$ is a Shintani lift of $\pi_v$.

The following theorem has already been proved for cuspidal representations by Arthur and Clozel ([AC]).\\
\ \\
{\bf Theorem E.}

{\it {\rm (a)} Let $\pi$ be an automorphic discrete series (resp. cuspidal) representation of $GL_n(\aa_{\F})$.   Set $r:=\frac{l}{m(\pi)}$. Then $n$ is divisible by $r$ and there exists a unique $\s^{r}$-stable discrete series (resp. cuspidal) representation $\pi'$ of $GL_{\frac{n}{r}}(\aa_{\E})$ such that 
$$\Pi=\pi'\tt  {\pi'}^\s \tt...\tt  {\pi'}^{\s^{r-1}}$$ 
is a $\s$-stable automorphic representation of $GL_n(\aa_{\E})$ and a Shintani lift of $\pi$. Moreover, $r(\pi')=r$.

{\rm (b)}  Let 
$$\Pi=\pi'\tt  {\pi'}^\s \tt...\tt  {\pi'}^{\s^{r-1}}$$ 
be a representation of  $GL_{n}(\aa_{\E})$, where $\pi'$ is a discrete series (resp. cuspidal) representation of  $GL_{\frac{n}{r}}(\aa_{\E})$ such that $r(\pi')=r$. Then $\Pi$ is the Shintani lift of some discrete series (resp. cuspidal) representation of $GL_n(\aa_{\F})$.}\\
\ \\
{\bf Remark.} The fibers of the lift can easily be described as in the local setting, following [AC] or [He1].

\section{Proofs}\label{proofs}

Until Section \ref{sepdiscser}, $E/F$ is an extension of local fields .

\subsection{Shintani lift for characters and twists}\label{liftcharacters}

Let $X(GL_n(F))$ (resp. $X(GL_n(E))$) be the group of smooth characters of $GL_n(F)$ (resp. of $GL_n(E)$). If $F$ is non-Archimedean, let $X^0(GL_n(F))$, $X^0(GL_n(E))$ the subgroups of unramified characters. Recall the characters $\chi\in X^0(GL_n(F))$ are of the form $g\mapsto |\det g|_F^z$, with $z$ a complex number determined by $\chi$ up to $\frac {2i\pi}{\ln q_F}\z$.

\begin{prop}\label{liftchar}
{\rm (a)} Every character  $\chi\in X(GL_n(F))$, has a Shintani lift  $\chi_{E}\in X(GL_n(E))$ defined by $\chi_{E}(g)=\chi(\nn g)$. The \nio\ is identity. 

{\rm (b)} If $F$ is non-Archimedean, then every unramified character $|\det|_F^z$, $z\in \cc$, lifts to $|\det|_E^z$. 
The lift is a surjective group morphism from $X^0(GL_n(F))$ to $X^0(GL_n(E))$. 

{\rm (c)} If $\chi\in X(GL_n(F))$ and $\pi$ is an irreducible representation of $GL_n(F)$ which has a Shintani lift $\Pi$ to $GL_n(E)$, then $\chi\pi$ has a  Shintani lift $\chi_{E}\Pi$.
\end{prop}
\ \\
{\bf Proof.} (a) In the {\bf Example} at the end of Section \ref{niosection} we explained why the \nio\ of $\chi_{E}$ is identity. One checks directly that the Shintani relation holds by the very definition of $\chi_{E}$.

(b) We have $|\det \nn(g)|_F=|N_{E/F}(\det(g))|_F=|\det (g)|_E$ since $|\ \ |_F\circ N_{E/F}=|\ \ |_E$. The rest is obvious. 

(c) is clear, we have already noted (end of Section \ref{niosection}) that the \nio\ is the same for $\Pi$ and $\chi_{E}\Pi$.\qed

\subsection{Shintani lift for products of representations}\label{liftproducts}

In the following we will study the lift of an irreducible product of representations (Proposition \ref{liftprod} below). 
It is more complicated than in the case where the representations are generic.

Let $\tau_i$, $i\in\k$, be an essentially tempered representation of $GL_{d_i}(E)$, $\tau_i=\nu^{e_i}\tau_i^u$ as in section \ref{notation}. Assume $\sum_{i=1}^k d_i=n$. Let $W$ be the space of the representation $\otimes_{i=1}^k\tau_i$ of the Levi subgroup $\prod_{i=1}^k GL_{d_i}(E)$ of $GL_n(E)$.
Set $D:=\tau_1\tt\tau_2\tt...\tt\tau_k$. 
Let $s\in\S_k$. Consider the representation $\tau_{s^{-1}(1)}\otimes\tau_{s^{-1}(2)}\otimes...\otimes\tau_{s^{-1}(k)}$ of the Levi subgroup $\prod_{i=1}^k GL_{d_{s^{-1}(i)}}(E)$ defined, to be coherent with the intertwining operators theory, on the same space $W$. Let $D_s$ be the representation $\tau_{s^{-1}(1)}\tt\tau_{s^{-1}(2)}\tt...\tt\tau_{s^{-1}(k)}$.
Then the normalized intertwining operator $N(s,D, (e_1,e_2,...,e_k)):D\to D_s$ from [MW1] page 607, is, by the property (2) there, well-defined and non-zero. Let us call it simply $N_s$ here -- its definition is recalled in the proof below.

\def\zz{{\underline{z}}}


\begin{lemme}\label{diagcommut}
(a) The diagram

$$
\begin{array}{ccccc}
 D & \xrightarrow{\s} &\s(D) \\
N_s\big\downarrow && \big\downarrow N_s \\
 D_s & \xrightarrow{\s} & \s(D_s) \\
\end{array}
$$
is commutative.

(b) Assume $\tau_i$ is $\s$-stable for $i=1,2,...,k$, and let $J:D\to \s(D)$ and $J_s:D_s\to \s(D_s)$ be the intertwining operators induced from the $I^{gen}_{\tau_i}$. Then the following diagram 
$$
\begin{array}{ccccc}
 D & \xrightarrow{J} &\s(D) \\
N_s\big\downarrow && \big\downarrow N_s \\
 D_s & \xrightarrow{J_s} & \s(D_s) \\
\end{array}
$$
is commutative.
\end{lemme}

\ \\
{\bf Proof.} (a) 
$N_s$ is the evaluation at $\zz=(e_1,e_2,...,e_k)$ of a product $a_sM_s$ where $M_s$ -- the unnormalized intertwining operator -- depends meromorphically on $\zz=(z_1,z_2,...,z_k)\in\cc^k$ and $a_s$  -- the normalizing factor -- is a meromorphic function of $\zz$. Precisely we have:
$$a_s(\zz)=\prod L(\tau_i^u\tt \check{\tau}_j^u, z_i-z_j)^{-1}
L(\tau_i^u\tt \check{\tau}_j^u, 1+z_i-z_j)
\varepsilon(\tau_i^u\tt \check{\tau}_j^u, z_i-z_j,\psi_E)$$
where the product is taken over inversions in $s$, that is over pairs of integers $(i,j)$, $1\leq i< j\leq k$ such that $s(i)>s(j)$. The $L-$ and $\varepsilon-$ factors are those defined by Shahidi [S1] or Jacquet, Piatetski-Shapiro and Shalika [JP-SS] (they are the same by [S2]). Here $\psi_E=\psi\circ\tr_{E/F}$ where $\psi$ is a non-trivial character of $F$. The $L-$ and $\varepsilon-$ factors for pairs depend only on the isomorphism classes of the representations involved -- indeed they are defined using their Whittaker models with respect to $\psi_E$ -- so the normalizing factor $a_s(\zz)$ stays the same if we replace each $\tau_i$ by $\s(\tau_i)$.\\

On the other hand the unnormalized intertwining operator $M_s(\zz)$ is given by an integral, which we now recall.

We let $L$ be the block-diagonal subgroup of $GL_n$ with consecutive blocks of size $d_1,d_2,...,d_k$ such that $\tau_i$ is a representation of $GL_{d_i}(E)$, $d_1+d_2+...+d_k=n$. Let $P$ be the upper block-triangular parabolic subgroup of $GL_n$ with Levi component $L$, and $U$ its unipotent radical. 

Gather the vectors in the canonical basis $(v_1,v_2,...,v_n)$ of $E^n$ in successive strings $S_1,S_2,...,S_k$, of size $d_1$, $d_2$,...,$d_k$ and let $w=w_s$ be the permutation $S_{s^{-1}(1)},S_{s^{-1}(2)},....,S_{s^{-1}(k)}$ of this basis (permute strings keeping the same order inside a given string); we also write $w$ for the corresponding permutation matrix of $GL_n(E)$. For example, when $k=2$ and $s=(1,2)$, $w$ is the matrix 
$\begin{pmatrix}
0 & I_{d_2}\\
I_{d_1} & 0\\
\end{pmatrix}
$.
Then $L'=wLw^{-1}$ is the block diagonal subgroup of $GL_n$ with consecutive sizes $d_{s^{-1}(1)}, d_{s^{-1}(2)},...,d_{s^{-1}(k)}$. Let $P'$ be the upper block-triangular subgroup with Levi component $L'$ and $U'$ its unipotent radical. Let $f$ be a function in the space\footnote{We should write $V_{\zz}$ where $\zz=(z_1,z_2,...,z_k)$, but, as is customary, we identify all these spaces by taking restriction to $K_E$. The operator $\s$ is compatible with this identification.} 
$V$ of $\nu^{z_1}\tau_1^u\tt\nu^{z_2}\tau_2^u\tt...\tt\nu^{z_k}\tau_k^u$. Consider the function $\phi: g\mapsto f(w^{-1}g)$ on $GL_n(E)$. 
For $l'\in L'=wLw^{-1}$ we have 
$$\phi(l'g)=f(w^{-1}l' g)=f(w^{-1}l' w w^{-1}g)=\phi_{\zz}(w^{-1}l'w)\phi(g)$$
where $\phi_\zz$ is the representation 
$\nu^{z_1}\tau_1^u\otimes \nu^{z_2}\tau_2^u\otimes ...\otimes \nu^{z_k}\tau_k^u$ of $GL_{d_1}(E)\tt GL_{d_2}(E)\tt ...\tt GL_{d_k}(E)=L$ 
so that $\phi(l'g)=\phi'_\zz (l')\phi(g)$ where $\phi'_\zz$ is the representation  
$\nu^{z_{s^{-1}(1)}}\tau_{s^{-1}(1)}^u\otimes \nu^{z_{s^{-1}(2)}}\tau_{s^{-1}(2)}^u\otimes ...\otimes \nu^{z_{s^{-1}(k)}}\tau_{s^{-1}(k)}^u$ 
of 
$GL_{d_{s^{-1}(1)}}(E)\tt GL_{d_{s^{-1}(2)}}(E)\tt ...\tt GL_{d_{s^{-1}(k)}}(E)=L'$. 
Moreover, $\phi$ is left invariant under $wUw^{-1}$, so that we can consider the integral
$$\int_{wUw^{-1}\cap U'\bc U'}\phi(u'g)du'$$ 
where $du'$ is an invariant measure on the quotient. 
It is proved by Shahidi [S2] that this integral converges when $\zz$ is in some open subset of $\cc^k$ and can be extended meromorphically to the whole of $\cc^k$, yielding an intertwining operator $M_s(\zz)$ from $\nu^{z_1}\tau_1^u\tt\nu^{z_2}\tau_2^u\tt...\tt\nu^{z_k}\tau_k^u$ 
to 
$\nu^{z_{s^{-1}(1)}}\tau_{s^{-1}(1)}^u\tt \nu^{z_{s^{-1}(2)}}\tau_{s^{-1}(2)}^u\tt ...\tt \nu^{z_{s^{-1}(k)}}\tau_{s^{-1}(k)}^u$. 
Since the Haar measure $du'$ is $\s$-invariant, the intertwining operator $M_s(\zz)$ commutes with $\s$. This proves (a).\\
\ \\

(b) Assume now each $\tau_i$ is $\s$-stable and comes with its \nio\ $I^{gen}_{\tau_i}:\tau_i\to \tau_i^\s$ ($\tau_i$ is generic). 
That gives a linear map $I:W\to W$, $I=\otimes_{i=1}^k I_{\tau_i}^{gen}$.
Then $J :D\to \s(D)$ and $J_s :D_s\to \s(D_s)$ are given by $J(f)(g)=I(f(g))$ and $J_s(h)(g)=I(h(g))$ for functions $f\in D$ and $h\in D_s$. 
We have to show that the diagram 

$$
\begin{array}{ccccc}
 D & \xrightarrow{J} &\s(D) \\
N_s\big\downarrow && \big\downarrow N_s \\
 D_s & \xrightarrow{J_s} & \s(D_s) \\
\end{array}
$$
is commutative.

We have to ``check the diagram for $\zz$ where integrals converge". Notice that $J$ and $J_s$ are compatible with the identification of spaces $V_{\zz}$ through restriction to $K_E$.
We have already noted that the operator $a_s(\zz)$ does not change when we replace the $\tau_i$ by $\s(\tau_i)$. 
Moreover, for $g\in GL_n(E)$ and $f\in D$, one has, for $\zz$ where integrals converge:
$$M_s(\zz)(J(f))= \int J (f)(w^{-1}u' g)du'=\int I((f(w^{-1}u'g))du'$$
$$=I(\int f(w^{-1}u'g)du')=J_s(M_s(\zz)(f))$$ 
(we used the fact that $I$ is linear and commutes with integrals). That shows that the diagram is commutative by unique meromorphic extension property.\qed
\ \\

For the definition of the p.i.m.o. in the next Lemma see Section \ref{niosection}.
\begin{prop}\label{LgdsQuot}
Let $\tau_i$, $1\leq i\leq k$, be essentially tempered representations of $GL_{n_i}(E)$ such that $\sum_i n_i=n$, and assume all $\tau_i$ are $\s$-stable.
Then the $\s$-operator on $\pi:=Lg(\tau_1,\tau_2,...\tau_k)$ obtained from the \nio\ of $\tau_1\otimes\tau_2\otimes ...\otimes\tau_k$ by p.i.m.o. is the \nio\ $I_\pi$ of $\pi$.
\end{prop}

\ \\
{\bf Proof.} Let us take, in the previous Lemma, $D_s$ to be $\tau_1\t\tau_2\t ...\t\tau_k$ in this order (and $D$ the product in standard order). Then the intertwining operator $N_s$ is known to contain $\pi$ in its image (theory of the Langlands quotient), with multiplicity one. It induces (as explained in the Appendix, in the paragraph just after the Remark) an intertwining operator $\tilde{N}_s$ of $\pi$.

By definition, the \nio\ $I_\pi$ of $\pi$ is the intertwining operator between $\pi$ and $\pi^\s$ obtained by restriction of $\s^{-1}\circ J$. 
The operator $\s^{-1}\circ J_s$ induces also by restriction an intertwining operator of $\pi$ with $\pi^\s$ (see the Appendix) which we denote $I(\pi)$. The preceding Lemma implies $\tilde{N}_s I_\pi \tilde{N}_s^{-1}=I(\pi)$, so $I(\pi)$ is the \nio\ of $\pi$ by Lemma \ref{nio}.\qed
\ \\

Let now $\Pi_1$ and $\Pi_2$ be two irreducible representations of $GL_{n_1}(E)$ and $GL_{n_2}(E)$, $n_1+n_2=n$, such that $\Pi:=\Pi_1\tt\Pi_2$ is irreducible.
Let $D_1=\tau_1\tt\tau_2\tt...\tt\tau_k$, $D_2=\tau'_1\tt\tau'_2\tt...\tt\tau'_{k'}$ be  standard representations such that $\Pi_1$ is the Langlands quotient of $D_1$ and $\Pi_2$ is the Langlands quotient of $D_2$ -- it is understood that 
the $\tau_i,\tau'_i$ are essentially tempered representations. 

Define $\tau''_i$, for $1\leq i\leq k+k'$ by $\tau''_i=\tau_i$ if $1\leq i\leq k$ and $\tau''_i=\tau'_{i-k}$ if $k+1\leq i\leq k+k'$.
Let $s\in \S_{k+k'}$ be such that $(\tau''_{s(i)})_{1\leq i\leq k+k'}$ is in standard order. 
Let $D$ be the standard representation which is the ordered product of the elements of $(\tau''_{s(i)})_{1\leq i\leq k+k'}$. By transitivity of the parabolic induction, we have $D_1\tt D_2=D_{s}$. In [Ta3], Prop. 2.2 and 2.3,  it is shown that, in general, $\Pi_1\t \Pi_2$ has always, even if not irreducible, an irreducible subquotient isomorphic to the Langlands quotient of $D$. 
In particular, here, $\Pi_1\t\Pi_2$ is isomorphic to the Langlands quotient of $D$. 

Assume that $\Pi_1$ and $\Pi_2$ are $\s$-stable. Then all the $\tau''_i$ are $\s$-stable. Proposition \ref{LgdsQuot} implies that the \nio\  of $\Pi_1\t\Pi_2$ is equal to the $\s$-operator on $\Pi_1\t\Pi_2$ obtained from the \nio\ $I_{\tau''_1}^{gen}\otimes I_{\tau''_2}^{gen}\otimes ...\otimes I_{\tau''_{k+k'}}^{gen}$ by p.i.m.o.. But then, by Proposition \ref{transitivite} (b), this is the $\s$-operator on $\Pi_1\t\Pi_2$ obtained from
$I_{\Pi_1}\otimes I_{\Pi_2}$ by the parabolic induction procedure.

We have proved the

\begin{prop}\label{nioprod}
Let $\Pi_1$ and $\Pi_2$ be two $\s$-stable irreducible representations of $GL_{n_1}(E)$ and $GL_{n_2}(E)$ such that $\Pi_1\tt\Pi_2$ is irreducible. Then $I_{\Pi_1\tt\Pi_2}$ is obtained by the parabolic induction procedure from $I_{\Pi_1}\otimes I_{\Pi_2}$.
\end{prop}

As a consequence, we get the following important proposition:

\begin{prop}\label{liftprod}
Let $\pi_1$ and $\pi_2$ be two irreducible representations of $GL_{n_1}(F)$ and $GL_{n_2}(F)$. Let $\Pi_1$ and $\Pi_2$ be irreducible $\s$-stable representations of $GL_n(E)$ such that $\Pi_i$ and $\pi_i$ verify the Shintani relation for $i=1,2$. Assume that $\pi_1\t\pi_2$ and  $\Pi_1\tt\Pi_2$ are irreducible. Then $\Pi_1\tt\Pi_2$ and $\pi_1\tt\pi_2$ verify the Shintani relation.
\end{prop}
\ \\
{\bf Proof.} By a result of Clozel ([Cl] Theorem 2), $\pi_1\tt\pi_2$ and $\Pi_1\tt\Pi_2$ verify the Shintani relation for the intertwining operator  $\s^{-1} (I_{\Pi_1}\tt I_{\Pi_2})$. We have just proved that this is also the \nio\ of $\Pi_1\tt\Pi_2$ (proposition \ref{nioprod}).\qed
\ \\
\ \\
{\bf Example.} Notice that Proposition \ref{liftprod} applies well to unitary representations, as their product is always irreducible. When the product of the lifts is not irreducible awkward things may happen. For example, consider the standard irreducible representation  $|\ \ |_F\tt \chi$ of $GL_2(F)$, where $\chi$ is the generator of $X(E/F)$ we fixed. The lift of $|\ \ |_F$ is $|\ \ |_E$, and the lift of $\chi$ is $1$. It turns out that $|\ \ |_E\tt 1$ is reducible (and each irreducible subquotient -- a character and an essentially square integrable representation -- is $\s$-stable and is a lift of an irreducible representation). The Langlands quotient $|\det|_E^{1/2}$ of $|\ \ |_E\tt 1$ is the base change of $|\ \ |_F\tt \chi$. But they do not verify the Shintani relation, because the character of $|\ \ |_F\tt \chi$, which is fully induced, vanishes on the set of elliptic elements, which always contains norms (every square of an element of $GL_2(F)$ is a norm) while the character $|\det|_E^{1/2}$ vanishes nowhere. (As we proved above, $|\det |_E^{1/2}$ is the Shintani lift of  $|\det|_F^{1/2}$.)\\
\ \\

The same proof of Proposition \ref{nioprod} implies, however, even when $\Pi_1\t\Pi_2$ is not irreducible:
\begin{prop}
The operator $I_{\Pi_1}\otimes I_{\Pi_2}$ induces by p.i.m.o. on the Langlands quotient of $D$, which appears in $\Pi_1\t\Pi_2$ with multiplicity one, the \nio.
\end{prop}
\ \\
{\bf Remark.}
Assume $\Pi_i$ is the base change of $\pi_i$, $i=1,2$ and $\pi_1\t\pi_2$ is irreducible. 
Then $D$ is the base change of $\pi_1\t \pi_2$. If, moreover, $\Pi_i$ is the Shintani lift of $\pi_i$ for $i=1,2$, 
 Proposition \ref{liftprod} shows that if $\Pi_1\t \Pi_2$ is irreducible then $D=\Pi_1\t\Pi_2$ is the Shintani lift of $\pi_1\t\pi_2$, while the {\bf Example} above shows that if $\Pi_1\t\Pi_2$ is reducible then $D$ may not be the Shintani lift of $\pi_1\t\pi_2$.

\subsection{Shintani lift for spherical unitary representations (Theorem B)} \label{sphericalunit} Here $F$ is a $p$-adic field. As a corollary of the lift of characters and of products, we get the Shintani lift for spherical unitary representations (proved in [AC] by other methods). Indeed, the classification of spherical representations combined with Tadi\'c's classification of unitary representations implies that $\pi$ is an irreducible spherical unitary representation  of $GL_n(F)$ if and only if $\pi$ is a product 
\begin{equation}\label{classifgenspher}
\pi=c_1\tt c'_1\tt c_2\tt c'_2\tt...\tt c_k\tt c'_k\tt d_1\tt d_2\tt ...\tt d_p,
\end{equation}
where 

- for each $i\in \{1,2,...,k\}$, $c_i=|\det|_F^{a+ib}$ and $c'_i=|\det|_F^{-a+ib}$ are unramified characters of some $GL_{n_i}(F)$ such that 
$a\in ]0,\frac{1}{2}[$, $b\in\r$, and 

- for each $j\in \{1,2,...,p\}$, $d_j$ is an unramified unitary character of some $GL_{m_j}(F)$, such that $2\sum_i n_i+\sum_j m_j=n$.

Moreover, the characters $c_i,c'_i$ and $d_j$ are determined by $\pi$ (up to permutation). The same is true for $GL_n(E)$. Gathering Proposition \ref{liftchar} and \ref{liftprod} we then get {\bf Theorem B}. The Shintani lift of $\pi$ is the corresponding product of characters of $GL_{n_i}(E),GL_{m_j}(E)$. Moreover, a spherical unitary representation $\pi'$ has the same lift as $\pi$ if and only if $\pi'=C_1\tt C'_1\tt C_2\tt C'_2\tt...\tt C_k\tt C'_k\tt D_1\tt D_2\tt ...\tt D_p$ where $C_i$ (resp. $C'_i$, $D_j$) are characters and have the same Shintani lift as $c_i$ (resp. $c'_i$, $d_j$). (We recall that two characters have the same lift if and only if they differ by a character in $X(E/F)$).\\
\ \\
{\bf Remark on the \nio.} The \nio\ of unramified characters is identity, and the maximal compact subgroup $K_E$ is stabilized by $\s$. It follows by (parabolic) induction that the \nio\  of $\pi$ induces by restriction the identity on the line of spherical vectors.\\ 
\ \\
{\bf Consequence.} We know that local components of automorphic cuspidal representations of $GL_n(\F)$, $\F$ global, are unitary spherical AND generic at almost every finite place of the global field. According to the classification of generic representations described in section \ref{notation}, all the $n_i$ and all the $m_j$ are then equal to $1$. Now every residual global representation $\pi$ is such that $\pi=MW(\rho,k)$ for $\rho$ some cuspidal representation. Hence the local component of $\pi$ at some place is $u(\tau,k)$, where $\tau$ is the local component of $\rho$ at that place (section \ref{classification}). So, at almost every place, the local component of $\pi$ (which is unitary and spherical) is a product like \ref{classifgenspher}, such that the $n_i$ and the $m_j$ are all equal to $k$. A consequence is that if the local component of $\pi$ at some finite place $v$ has the same Shintani lift as the local component at $v$ of some other residual representation $\pi'$, we have $\pi'=MW(\rho',k)$ where the local component of $\rho'$ at $v$ has the same lift as the local component of $\rho$ at $v$. This fact will be used later.\\
\ \\
{\bf Remark on the non-unitary case.} When $\pi$ is spherical but non-unitary, it is still an irreducible product $\chi_1\t \chi_2\t ...\t \chi_k$ where $\chi_i$, $1\leq i\leq k$, is an unramified character of some $GL_{n_i}(F)$. If the product $\chi_{1,E}\t \chi_{2,E}\t ...\t \chi_{k,E}$ of the Shintani lifts is  irreducible, then the same proof gives the Shintani lift of $\pi$. The Example following Proposition \ref{liftprod} shows the problems arising when the product of the Shintani lifts is reducible. One may easily verify if the conditions for the lift are fulfilled or not, it is just a question of linked Zelevinsky segments. For example, if we say that two lines of characters of $GL_n(F)$ 
$\{\nu^r\theta, r\in \r\}$ and $\{\nu^r\theta', r\in \r\}$, $\theta, \theta'$ unitary characters of $GL_1(F)$, are {\bf equivalent for base change} if there exist $\chi\in X(E/F)$ such that $\theta'=\chi\theta$, then if the cuspidal support of $\pi$ is included in a union of non-equivalent lines, the product of the Shintani lifts is irreducible and the base change is indeed a Shintani lift.

\def\res{{\rm res}}
\def\ind{{\rm ind}}

\subsection{Shintani lift for elliptic representations (Theorem C)}\label{sectelliptic}

In this section we assume that {\bf Theorem A} (b) is proved. {\bf Theorem A} is proved in [AC] and requires global methods. We shall come back to it in section \ref{trform}. We assume first $F$ is $p$-adic, then we treat the case $F\simeq \r$ and $E\simeq \cc$.

Let $\tau$ be an essentially square integrable representation of $GL_n(E)$. As recalled in the chapter \ref{notation} there exist a positive integer $k$, $k|n$, and a cuspidal representation $\rho$ of $GL_{\frac{n}{k}}(E)$ such that $\tau$ is  the unique irreducible subrepresentation of
$\ind_{L_0}^{GL_n(E)}S$
where $L_0$ is the standard Levi subgroup of $GL_n(E)$ with $k$ blocks of equal size $m:=\frac{n}{k}$ and 
$S=\nu^{k-1}\rho\otimes \nu^{k-2}\rho\otimes ...\otimes\rho$. 
The standard Levi subgroups of $GL_n(E)$ containing $L_0$  are parametrized by $I:={\mathcal P}(K)$ (set of subsets of $K$) where $K:=\{1,2,...,k-1\}$: 
if $L$ is $GL_{n_1 m}(E)\t GL_{n_2 m}(E)\t ... \t GL_{n_t m}(E)$, 
then we set $L=L_i$, where $i\in {\mathcal{P}}(K)$ is the complementary set of $\{n_1,n_1+n_2,n_1+n_2+n_3,...,n_1+n_2+...+n_{t-1}\}$ in $K$.  
Then $L_0=L_\emptyset$, $L_K=GL_n(E)$ and $L_i\subset L_j$ if and only if $i\subset j$. For all $i\in {\mathcal{P}}(K)$, let $\tau_i$ be the unique  irreducible subrepresentation of $\ind_{L_\emptyset}^{L_i} S$. For example $\tau_\emptyset=S$, $\tau_{K}=\tau$. 
Then $\tau_i$ is an essentially square integrable representation of $L_i$. 
Now let $X_{\emptyset}:=\ind_{L_\emptyset}^{L_K}S$ and, for all $i\in {\mathcal{P}}(K)$, $X_i$ the subrepresentation $\ind_{L_i}^{L_K}\tau_i$ of $X_{\emptyset}$. 
If $i,j\in I$ and $i\subset j$, then $X_j$ is a subrepresentation of $X_i$. 
Each $X_i$ has a unique irreducible quotient $\pi_i$, the Langlands quotient, denoted here $Lg(X_i)$.
It is known ([Ze], [BW] X 4.6) that the induced representation $X_{\emptyset}$ has exactly $2^{k-1}$ irreducible subquotients, the $\pi_i$ for $i\in {\mathcal{P}}(K)$, which appear with multiplicity one.
Moreover, $\pi_j$ is a subquotient of $X_i$ if and only if $i\subset j$.  
The elliptic representations of $GL_n(E)$ are exactly the representations $\pi_i$ constructed from essentially square integrable 
representations $\tau$ in this way ([Ba2] 2.5). These facts are true also for $GL_n(F)$ and on this group we will use the notations with a tilde ($\tilde{\pi}_i$ for example). 

\def\ttt{{\tilde{\tau}}}
\def\tpi{{\tilde{\pi}}}

Let $\ttt$ be an essentially square integrable representation of $GL_n(F)$ and set $r=l/m(\ttt)$. Let $u:=\tau\t\tau^\s\t\tau^{\s^2}\t...\t\tau^{\s^{r-1}}$ be the Shintani lift of $\ttt$ to $GL_n(E)$, where $\tau$ is an essentially square integrable representation of $GL_{n/r}(E)$ ({\bf Theorem A}).

\begin{prop}
Let $i\in {\mathcal{P}}(K)$, and $\tpi_i$ be the elliptic representation of $GL_n(F)$ associated to $\ttt$ and $i$ as before. Then $\tpi_i$ has a Shintani lift. Its lift is the  representation $\pi_{(i)}=\pi_i\t\pi_i^\s\t\pi_i^{\s^2}\t...\t\pi_i^{\s^{r-1}}$ of $GL_n(E)$, where $\pi_i$ is the elliptic representation associated to $\tau$ and $i$ as before.

\end{prop}
\ \\
{\bf Proof.}  Let us first treat the case when $m(\ttt)=l$, i.e. $r=1$ and $\tau$ is $\s$-stable. We use the notation above. Because $\tau_i$ is the Shintani lift of $\tilde{\tau}_i$, $\pi_{(i)}=\pi_i$ is the base change of $\tilde{\pi}_i$.
Also, $S$ is $\s$-stable and generic, and has a \nio\ $I_S^{gen}$. 
Let $i\in {\mathcal{P}}(K)$. The $\s$-operator  
 obtained from $I_S^{gen}$ by p.i.m.o. on $\tau_i$ is the \nio\ $I_{\tau_i}^{gen}$ of $\tau_i$, by Proposition \ref{transitivity}.
 The operator $I_{\tau_i}^{gen}$ induces an intertwining operator $I_{X_i}$ on $X_i$ by the parabolic induction procedure.
 By definition, the \nio\ $I_{\pi_i}$ of $\pi_i$ is obtained from $I_{X_i}$ by the multiplicity one property.
Then the $\s$-operator  obtained from $I_S^{gen}$ by p.i.m.o. on $\pi_i$ is the \nio\ $I_{\pi_i}$ of $\pi_i$, by Proposition \ref{transitivite} (b).

Because $\tau_i$ is the Shintani lift of $\ttt_i$, we have $\tr  I_{X_i} X_i(f')=\tr\tilde{X}_i(f)$ for functions $f\lra f'$ (by Theorem 2 in [Cl]).
Fix functions $f,f'$ such that $f\lra f'$.  
We want to show that, for all $i\in {\mathcal{P}}(K)$, 
$$\tr I_{\pi_i}\pi_i(f')=\tr\tilde{\pi}_i(f).$$ 
That follows by decreasing induction on $i$ from the formula:
$$\sum_{i\subset j}\tr \tilde{\pi}_j (f) = \tr \tilde{X}_i (f) = \tr  I_{X_i} X_i(f') = \sum_{i\subset j} \tr  I_{\pi_j}\pi_j(f').$$
\ \\ 

Let us move now to the general case. Let $\Theta$ be the induced representation 
$$\nu^{k-1}\rho\t\nu^{k-1}\rho^\s\t ...\t\nu^{k-1}\rho^{\s^{r-1}}\ \t\ 
\nu^{k-2} \rho\t\nu^{k-2}\rho^\s\t\  ...\t\nu^{k-2}\rho^{\s^{r-1}}\t...\ \t\ 
\rho\t\rho^\s\t ...\t\rho^{\s^{r-1}}\simeq$$
$$\simeq \nu^{k-1}u\t\nu^{k-2} u\t ...\t u,$$
where $u:=\rho\t\rho^\s\t ...\t\rho^{\s^{r-1}}$ is obviously essentially tempered.

By Proposition 8.5 of [Ze], the representation $\nu^{a}\rho^{\s^j}\t\nu^{b}\rho^{\s^{j'}}$ is irreducible and isomorphic to
$\nu^{b}\rho^{\s^{j'}}\t \nu^{a}\rho^{\s^j}$
for any $0\leq j<j'\leq r-1$ and any $a,b\in\r$ (because $\rho^{\s^j}$ and $\rho^{\s^{j'}}$ are not isomorphic, of the same exponent, so they are on different lines).
So  $\Theta$ is isomorphic to

$$(\nu^{k-1}\rho\t\nu^{k-2}\rho\t ...\t\rho)\ \t\ (\nu^{k-1}\rho^\s\t\nu^{k-2}\rho^\s\t ...\t \rho^\s)\t ...\t(\nu^{k-1}\rho^{\s^{r-1}}\t \nu^{k-2}\rho^{\s^{r-1}}\t ...\t\rho^{\s^{r-1}}).$$

This implies that:

- the irreducible subquotients of $\Theta$ are of multiplicity one, of the form $\pi_{i_1}\t\pi_{i_2}^\s\t\pi_{i_3}^{\s^2}\t ...\t\pi_{i_r}^{\s^{r-1}}$, where $\pi_{i_j}$ are chosen among the $\pi_i$, $i\in \mathcal{P}(K)$,

- for each $i\in \mathcal{P}(K)$ we have a subrepresentation $X_{(i)}\simeq X_i\t X_i^\s\t ...\t X_i^{\s^{r-1}}$ of $\Theta$.

Let $I_\Theta$ be the $\s$-operator on $\Theta$ obtained from $I^{gen}_{\nu^{k-1} u}\otimes I^{gen}_{\nu^{k-2 } u}\otimes ...\otimes I^{gen}_u$ by the parabolic induction procedure. By transitivity of the induction functor, the irreducible subquotients of $X_{(i)}$ are the $\pi_{i_1}\t\pi_{i_2}^\s\t\pi_{i_3}^{\s^2}\t ...\t\pi_{i_r}^{\s^{r-1}}$ with $i\subset i_j$ for all $j$, $1\leq j\leq r$. In particular, any subrepresentation of $\Theta$ isomorphic to $X_{(i)}$ is equal to $X_{(i)}$ (see [Ze], Chapter 2, where it is proven that, for a representation with multiplicity free cuspidal support -- like $\Theta$ -- a submodule is determined by the set of isomorphism classes of its irreducible subquotients). So, as $X_{(i)}$ is $\s$-stable, $X_{(i)}$ is stable by $I_\Theta$.

For $i\in \mathcal{P}(K)$, set $\pi_{(i)}:=\pi_{i}\t\pi_{i}^\s\t\pi_{i}^{\s^2}\t ...\t\pi_{i}^{\s^{r-1}}$. The $\pi_{(i)}$, $i\in \mathcal{P}(K)$, are the irreducible subquotients of $\Theta$ which are $\s$-stable.
 Let    $I_\Theta (\pi_{(i)})$ be the $\s$-operator on $\pi_{(i)}$ obtained from $I_{X_{(i)}}$ by the multiplicity one property. We have the

\begin{lemme} 
$I_\Theta (\pi_{(i)})$ is the \nio\ $I_{\pi_{(i)}}$.
\end{lemme}
\ \\
{\bf Proof.} Recall that $\tau_i$ is the unique subrepresentation of 
$\nu^{k-1}\rho\t\nu^{k-2}\rho\t ...\t\rho$, that $\tau_i$ is essentially square integrable and $X_i=\ind_{L_i}^{L_K} \tau_i$. Then $\pi_i=Lg(X_i)$. 
Write $\tau_i=\tau_i^1\otimes \tau_i^2\otimes ...\otimes \tau_i^{m(i)}$, where $m(i)$ is the number of blocks of $L_i$ (not to be confused with $m(\rho)$ where $\rho$ is a representation) and the representations $\tau_i^1, \tau_i^2, ...,\tau_i^{m(i)}$ are essentially square integrable, in standard order. For $1\leq j\leq m(i)$, set 
$\tau_{(i)}^j:=\tau_i^j\t (\tau_i^j)^\s\t ...\t (\tau_i^j)^{\s^{r-1}}$. This is an essentially tempered representation, the $\tau_{(i)}^j$ are in standard order (decreasing with $j$), and $\pi_{(i)}$ is $Lg(\tau_{(i)}^1,\tau_{(i)}^2,...,\tau_{(i)}^{m(i)})$ (this is a consequence of [Ta3], Prop. 2.2 and 2.3, because $\pi_{(i)}$ is the irreducible product of the $\pi_i^{\s^t}$, $0\leq t\leq r-1$, and $\pi_i^{\s^t}=Lg(X_i^{\s^t}$)). In particular, $\pi_{(i)}$ is the base change of $\tilde{\pi}_i$.
The \nio\ $I_{\pi_{(i)}}$ is, by definition, obtained from $I_{\tau_{(i)}^1}\otimes I_{\tau_{(i)}^2}\otimes...\otimes   I_{\tau_{(i)}^{m(i)}}$ by p.i.m.o..

Let ${\alpha_i^j}$ be the length of the segment of $\tau_i^j$. Then $\nu^{k-1} u\t \nu^{k-2}u\ ...\t \nu^{k-\alpha_i^1}u$ is a subrepresentation of a representation induced from a segment of length $\alpha_i^1 r$ and has $\tau_{(i)}^1$ as a subquotient of multiplicity one. As $\tau_{(i)}^1$ is a generic representation, Proposition \ref{transitivity} shows that its \nio\
is obtained from  $I^{gen}_{\nu^{k-1} u}\otimes I^{gen}_{\nu^{k-2 } u}\otimes ...\otimes I^{gen}_{\nu^{k-\alpha_i^1} u}$ by p.i.m.o.. The same is true for the other $\tau_{(i)}^j$, $2\leq j\leq m(i)$. The lemma follows then by transitivity -- Proposition \ref{transitivite} (b).\qed
\ \\

Fix functions $f,f'$ such that $f\lra f'$. 
We claim that we have: 

$$\tr  I_{X_{(i)}} X_{(i)}(f') = \sum_{i\subset j} \tr  I_{\pi_{(j)}}\pi_{(j)}(f').$$ 

Indeed, let $0\subset \pi_{K} \subset U_2\subset ...\subset U_m=X_{(i)}$  a Jordan-H\"{o}lder series for the action of $GL_n(E)$ via $\Theta$ {\it and} $I_\Theta$ (obviously of finite length), meaning that all the submodules in the series are stable both by $\Theta$ and $I_\Theta$ and the consecutive quotients are irreducible for this action (see [Bou] I.1.4.7 for the general version of Jordan-H\"{o}lder theory for groups with operators). On the one hand, the trace of $I_{X_{(i)}} X_{(i)}(f')$ is the sum of the trace on the quotients. On the other hand, by the uniqueness of the composition series (see [Bou] p. 43, Theorem 6) we have that all the $\pi_{(j)}$, $i\subset j$, being stable by $I_\Theta$ and of multiplicity one, appear with multiplicity one in the composition series. Moreover the trace is null on the other quotients : consider an irreducible subrepresentation $\varepsilon$ (for the action of $GL_n(E)$) of this quotient $U_{l+1}/U_l$. Then $\varepsilon$ is isomorphic to some $\pi_{i_1}\t\pi_{i_2}^\s\t\pi_{i_3}^{\s^2}\t ...\t\pi_{i_r}^{\s^{r-1}}$, with the index $i_j$ not all equal. Then $I_\Theta$ sends $\varepsilon$ to some {\it different} irreducible subrepresentation of $X_{(i)}$. The quotient $U_{l+1}/U_l$ is the sum of the conjugates of $\varepsilon$ under $I_\Theta$, and if there is more than one such conjugate and they are permuted by $I_\Theta$ without fixed point, the trace is null.

Now the proof goes as in the case $m(\tau)=1$: we have on the group $GL_n(F)$ the relation
$$\tr \tilde{X}_i(f)=\sum_{i\subset j}\tr \tilde{\pi}_j (f),$$
and by the same arguments we have  $\tr \tilde{X}_i(f)=\tr  I_{X_{(i)}} X_{(i)}(f')$. We conclude by decreasing induction on $i$.\qed
\ \\
\ \\

A few words about claim (b) of {\bf Theorem C}, which is now easy to prove. Let $\xi_E$ be as in claim (b); using the preceding construction -- and notation -- we may assume that $\xi'=:\pi_i$ is the Langlands quotient of $X_i$ as in the preceding proof. Then $\tau_{(i)}$ is the Shintani lift of some $\tilde{\tau_i}$ by the base change theory in the tempered case, and $\tilde{\pi}_i$ is an elliptic representation of $GL_n(F)$. The preceding Proposition prove that the lift of $\tilde{\pi}_i$ is $\Xi=\xi'\t{\xi'}^\s\t ...{\xi'}^{\s^{r-1}}$.\\

\ \\
{\bf Remark.} The Shintani lift of an elliptic representation is obtained here by local methods using the Shintani lift of square integrable representations obtained by global methods. Notice that we know (by local methods) the Shintani lift of characters. This may serve as a first step of a proof by induction of the Shintani lift of  the Steinberg representation and its twists by characters (by purely local methods). The proof would be similar to the previous one. One has to take $\rho$ a character of $GL_1(F)$ (then $m(\rho)=1$), reverse the order ($X_{\emptyset}=\rho\t\nu\rho\t ...\t\nu^{k-1}\rho$), invoke Proposition \ref{LgdsQuot} to show that the induced intertwining operators on the irreducible subquotients are their \nio s and apply the induction for increasing $i$.\\
\ \\
{\bf The Archimedean case.} The group $GL_n(\cc)$ has elliptic elements if and only if $n=1$. In this case all the irreducible smooth representations are characters.
They are of the form $\xi_{k,\alpha}(z)=|z|^\alpha z^k$ for $\a\in \cc$ and $k\in\z$. Such a character is stable by conjugation if and only if $k=0$. We then denote it simply $\xi_\a$.

\def\b{\beta}
\def\sgn{{\rm sgn}}
The group $GL_n(\r)$ has elliptic elements and elliptic representations if and only if $n=1$ or $n=2$ (because irreducible polynomials of $\r[X]$ are of degree $1$ or $2$). For $n=1$, the elliptic representations are characters. They are of the form $\chi_\a(x)=|x|^\a$ or $\chi'_\a(x)=|x|^\a\sgn$, $\a\in\cc$, and $\sgn$ is the sign character. It is known that 

(a)  induced representations $\chi_\a\t\chi'_\b$ 
and $\chi'_\a\t \chi_\b$ are reducible if and only if $\a-\b$ is an even integer,

(b)  induced representations $\chi_\a\t\chi_\b$ and $\chi'_\a\t \chi'_\b$ are reducible if and only if $\a-\b$ is an odd integer. 

Let $m\in \n^*$. 

If $m$ is even, $\chi_{m/2}\t \chi'_{-m/2}$ has two irreducible subquotients, its Langlands quotient $r(m)$ and an essentially square integrable representation $\d(m)$. The representation $\d(m)$ is stable by multiplication with the character sign, while $r(m)$ is not. So, after semisimplification, one gets in the Grothendieck group: $\chi'_{m/2}\t\chi_{-m/2}=\d(m)+r'(m)$, $m\in\n^*$, where $r'(m)=\sgn \otimes r(m)$.

When $m$ is odd, a completely parallel situation occurs : $\chi_{m/2}\t \chi_{-m/2}$ has two irreducible subquotients, its Langlands quotient $r(m)$ and an essentially square integrable representation $\d(m)$. The representation $\d(m)$ is stable by multiplication with the character sign, while $r(m)$ is not. So, after semisimplification, one gets in the Grothendieck group: $\chi'_{m/2}\t\chi'_{-m/2}=\d(m)+r'(m)$, $m\in\n^*$, where $r'(m)=\sgn \otimes r(m)$.

A character of $GL_2(\r)$ is of the form $\chi_\a\circ \det$ or $\chi'_\a\circ \det$, $\a\in\cc$. The representations $\d(m),r(m)$, $m\in\n^*$, are elliptic. Every elliptic representation of $GL_2(\r)$ is (up to isomorphism) a twist of some $\d(m)$ or $r(m)$ with a character. Indeed, all the other irreducible representations are fully induced from characters of the diagonal torus. It is known that $r(m)$ is of dimension $m$, and $r(1)$ is the trivial representation. 

The only interesting cases of Archimedean Shintani lift for elliptic representations are then when $E/F$ is (up to isomorphism) $\cc/\r$ and $n=1$ or $n=2$. In both cases, the sign characters  (i.e. $\sgn$ when $n=1$ and $\sgn\circ\det$ when $n=2$) are trivial on the norms. In the case $n=1$, the Shintani lifts of $\chi_\a$ and $\chi'_\a$ are the same, equal to $\xi_{2\a}$, and in the case $n=2$, the Shintani lifts of $\chi_\a\circ\det$ and $\chi'_\a\circ\det$ are the same, equal to $\xi_{2\a}\circ\det$. This is straightforward checking. 

So, in the case $n=2$, it is enough to explain Shintani lift for representations $\d(m)$ and $r(m)$, $m\in\n^*$, as the lifts of the other elliptic representations of $GL_2(\r)$ are then obtained by twists with characters. The induced representation $\xi_{m}\t \xi_{-m}$ to $GL_2(\cc)$ has two irreducible subquotients: the Langlands quotient $R(m)$ and the irreducible subrepresentation $\Delta(m):=\theta\t\bar{\theta}$, where $\theta(z):=\frac{z^m}{|z|^m}$ (i.e. $\theta=\xi_{-m,m}$) and $\bar{\theta}$ is its complex conjugate (notice that both characters are unitary). Since the \nio s of the characters $\xi_m$ and $\xi_{-m}$ are trivial, they induce (see the construction in Section \ref{niosection}) the operator $f\mapsto \tilde{f}$ on the space of functions of the induced representation, where $\tilde{f}(z)=f(\bar{z})$ for all $z\in GL_2(\cc)$. This operator 
induces on the Langlands quotient $R(m)$, by definition, the \nio\ of $R(m)$, and on $\Delta(m)$, because it is generic -- Proposition \ref{transitivity} --, the  \nio\ of $\D(m)$. By [AC] we know that the Shintani lift of $\d(m)$ is $\Delta(m)$. 
Then the same method as for $p$-adic groups (this section) shows that the Shintani lift of $r(m)$ is $R(m)$. See also [La], chapter 7.

\subsection{Separating discrete series}\label{sepdiscser}

\def\ind{{\rm{ind}}}
This section is global.
If $k|n$, we denote $L_k$ the Levi subgroup of $GL_n(\aa_\F)$ of diagonal matrices by $k$ equal blocks. If $\d$ is an automorphic discrete series of $GL_\frac{n}{k}(\aa_\F)$, we denote $\d^k$ the automorphic representation obtained by parabolic induction from $L_k$ to $GL_n(\aa_\F)$ (with respect to the upper triangular parabolic subgroup as in the local setting) of the tensor product of $k$ copies of $\d$. As local components of $\d$ are unitary, the induced representation is irreducible, locally and hence globally.\\
\ \\
\begin{prop}\label{separating}
Let $\pi$, $\pi'$ be two representations of $GL_n(\aa_{\F})$ such that $\pi=\d^m$  and $\pi'=\d'^{m'}$ for automorphic discrete series $\d$ of $GL_{\frac{n}{m}}(\aa_{\F})$ and $\d'$ of $GL_{\frac{n}{m'}}(\aa_{\F})$. Let $S$ be a finite set of places of $\F$ containing all the infinite places and all the finite places $v$ such that either $\E/\F$ is ramified, or $\pi_v$ or $\pi'_v$ is not a spherical representation. Assume that for every place $v\notin S$ one has $\tr\pi_v(b(f))=\tr\pi'_v(b(f))$ for all $f\in H^0(G_{\E_w})$, where $w$ is a place of $\E$ such that $w|v$. Then $m=m'$ and there exists a character $\chi\in X(\E/\F)$ such that $\d=\chi\d'$ and $\pi'=\chi \pi$.
\end{prop}

\ \\
{\bf Proof.} 
Let us write $\d=MW(\rho,k)$ and $\d'=MW(\rho',k')$ with $\rho,\rho'$ cuspidal. We will show that $k=k'$ and $\rho_v^m$ and ${\rho'_v}^{m'}$ have the same base change for every $v\notin S$. All the arguments are in the {\bf Consequence}, section \ref{sphericalunit}. Indeed, $\d_v^m$ is a product of unitary unramified characters of $GL_k(\E_v)$, while ${\d'_v}^m$ is a product of unitary unramified characters of $GL_{k'}(\E_v)$, hence $k=k'$ (by unicity in the classification of unitary representations). As explained in the {\bf Consequence}, the automorphic cuspidal representations $\rho^m$ and $\rho'^{m'}$ have the same lift at almost every place.

We now reason as in [AC] Theorem 3.1, page 201, where the case $m=m'=1$ is treated. Keeping their notation $\eta$ and $L^S$, we find, under our hypothesis, that  $\prod_{i=1}^l L^S(s,\rho \otimes \tilde{\rho}\otimes \eta^i)^m$ 
and  
$\prod_{i=1}^l L^S(s,\rho' \otimes \tilde{\rho'}\otimes \eta^i)^{m'}$
are equal.
Since the left-hand side has a pole at $s=1$ ([JS2]), so has the right-hand side, which implies
that $\rho'=\eta^i\rho$ for some $i\in \{1,2,...,l\}$; we then get $m=m'$ by looking at the order of the pole at $s=1$.\qed
\ \\
\ \\

\subsection{The twisted trace formula comparison}\label{trform}

The equality \ref{traceformula} below is the main theorem of [AC]. The proof occupies the first 200 pages of the book and is based on previous work of Arthur. 

Let $\F_\infty:= \F\otimes_{\q} \r$ be the product $\prod_{v}\F_v$ where $v$ runs over the infinite places of $\F$. Let $\mu$ be a unitary character
of $\F^\t_\infty$. We use the embedding of $\F^\t_\infty$ in
$\aa^\times_\F$ at infinite places to realize it as a subgroup of
the center $Z(\aa_{\F})$. 

Let $\mu_\E$ be the lift of $\mu$ to
$\E^\t_\infty=\prod_w \E^\t_w$ where $w$ runs over the set of infinite
places of $\E$, i.e. $\mu_E=\mu\circ \prod_{v {\text{ infinite}}}\prod_{w|v}{\bf N}_{\E_w/\F_v}$. We use the embedding of $\E^\t_\infty$
in $\aa^\times_\E$ at infinite places to realize it as a
subgroup of the center $Z(\aa_{\E})$.

Let ${\mathcal L}(\F)$ (resp. ${\mathcal L}(\E)$) be the set of $\F$-Levi subgroups of $GL_n(\F)$ 
(resp. $\E$-Levi subgroups of $GL_n(\E)$) containing the maximal diagonal torus.

The following formula \ref{traceformula} is the formula (4.1)=(4.2), page 203, in [AC], with $\s^{-1}$ in 
place\footnote{As $\s$ is any generator of $Gal(\E/\F)$ in [AC], one may switch to $\s^{-1}$ and the formula is still valid. In [AC] it is not clearly stated form which space to which space do the operators $\s$ and  $M(s,0)$ map. It is just a matter of convention, but the conditions $s\sigma\pi_\E=\pi_\E$ (page 207) and $s\pi_\E=\s \pi_\E$ (page 213), necessary for the term associated to $(s,\s)$ in the formula not to vanish, seem to correspond to conventions  opposite to one another. Here we work with our convention for $\s$ and for $M(s,0)$ (which is the one of [MW]).} 
of $\s$:

\begin{equation}\label{traceformula}
\sum_{L\in {\mathcal L}(\F)} |W_0^L||W_0^{GL_{n}(\F)}|^{-1}
\sum_{s\in W({\mathfrak a}_L)_{reg}} |\det(s-1)_{{\mathfrak a}_L^{GL_{n}(\F)} }|^{-1}
\tr(M_{L}^{GL_{n}(\aa_{\F})}(s,0)\ \circ\ \rho_{L,t,\mu}(f))=$$
$$l\, \sum_{L'\in {\mathcal L}(\E)} |W_0^{L'}||W_0^{GL_{n}(\E)}|^{-1}
\sum_{s\in W({\mathfrak a}_{L'})_{reg}} |\det(s-1)_{{\mathfrak a}_{L'}^{GL_{n}(\E)} }|^{-1}
\tr(M_{L'}^{GL_{n}(\aa_{\E})}(s,0)\ \circ\ \s^{-1}\ \circ \rho_{L',t,\mu_\E}(\phi))
\end{equation}
 where (see also [AC], page 132):

- $t\in \r_+$;

- $|W_0^L|$ is the cardinality of the Weyl group of $L$;

- ${\mathfrak a}_L$ is the real space $Hom(X(L)_\F,\r)$ where $X(L)_\F$ is
the lattice of rational characters of $L$; $W({\mathfrak a}_L)$ is the
Weyl group of ${\mathfrak a}_L$; ${\mathfrak a}_L^{GL_{n}}$ is the
quotient of ${\mathfrak a}_L$ by ${\mathfrak a}_{GL_{n}}$; 
$W({\mathfrak a}_L)_{reg}$ is the set of $s\in W({\mathfrak a}_L)$ such that
$\det(s-1)_{{\mathfrak a}_L^{GL_{n}}}\neq 0$;

- $\rho_{L,t,\mu}$ is the induced representation with respect to any
parabolic subgroup with Levi factor $L(\aa_\F)$ from the direct sum of
discrete series $\pi$ of $L(\aa_\F)$ such that $\pi$ is $\mu$-equivariant 
(i.e. the restriction of its central character to $\F_\infty^\t$ equals $\mu$)
and the imaginary part of the Archimedean infinitesimal character of
$\pi$ has norm $t$ ([AC], page 131-132); $\rho_{L',t,\mu_\E}$ is the corresponding representation when the field is $\E$.

- $M_L^{GL_n}(s,0)$ is the global intertwining operator associated to $s$ at the
point $0$; we sometimes denote it $M(s,0)$ when it is not necessary to specify the Levi subgroup;

- the operator $\s^{-1}$ is the operator $f\mapsto f\circ \s^{-1}$ in the space of $\rho_{L',t,\mu_\E}$ which is a space of functions stable by $\s$ (because of the choice of the central character). 

- $\phi$ and $f$ are associated. This means that $\phi=\otimes_v \phi_v$ and $f=\otimes_v f_v$, where $v$ runs over the places of $\F$, are such that for almost all $v$ where $\E_v/\F_v$ is unramified $\phi_v$ and $f_v$ are spherical and $b(\phi_v)=f_v$, and for the other places $v$, $\phi_v$ and $f_v$ match (chapter \ref{notglobal}).

Recall, when $k\in\n^*$ and $k|n$, $L_k\in {\mathcal L}(\F)$ and $L'_k\in {\mathcal L}(\E)$ are the Levi groups of diagonal matrices by $k$-blocks of the same size $\frac{n}{k}$. The set $W({\mathfrak a}_L)_{reg}$ (resp. $W({\mathfrak a}_{L'})_{reg}$) in the formula is empty, unless $L$ is conjugate to $L_k$ (resp. $L'$ is conjugate to $L'_k$) for some $k$ dividing $n$, and $s$ is a cycle of length $k$ permuting the blocks of $L_k$ (resp. $L'_k$). Moreover, it is shown in [AC], pages 207 to 209, that different cycles of length $k$ give the same contribution, and conjugate Levi subgroups give the same contribution to the trace formula \ref{traceformula}. So we will compute -- both left and right -- the contribution of the Levi subgroup $L_k$ (resp $L'_k$) and the 
cycle $s_k:=(k,k-1,k-2,...,1)$ then count the number of terms associated to that 
contribution\footnote{Once $s$ is fixed, only some terms, respecting symmetries, do not vanish in the trace formula. With our choice of cycle $s:=(k,k-1,k-2,...,1)$, these terms are of type $D:=\Delta\t\Delta^\s\times \Delta ^{\s^2}\times ...\times \Delta ^{\s^{k-1}}$, as in [AC].} .

If $\d_i$, $1\leq i\leq k$, are automorphic discrete series of $GL_{\frac{n}{k}}(\aa_\F)$, then we let $D:=\d_1\t\d_2\t ...\t\d_k$ be the automorphic representation of $GL_n(\aa_\F)$ parabolically induced from $\d_1\otimes\d_2\otimes ...\otimes\d_k$ (this induced representation is irreducible, because the local component at any place of $\d_i$, $1\leq i\leq k$, is unitary).
Then the representation $\rho_{L_k,t,\mu}$ is the sum of representations of the type of $D$.  
For our choice of $s:=s_k$, the intertwining operator $M(s,0)$ intertwines the space of $D:=\d_1\t\d_2\t ...\t\d_k$ with the space of 
$D_s:=\d_{s^{-1}(1)}\t\d_{s^{-1}(2)}\t\d_{s^{-1}(k)}=
\d_2\t\d_3\t ...\t\d_{k}\t\d_1$. These two spaces are either equal -- if and only if $\d_1=\d_2=...=\d_k$ -- or, else, disjoint. So the trace of $M(s_k,0)$ is zero unless the representation $D$ is of type $\d^k$.
 
The representation $\rho_{L'_k,t,\mu_\E}$ is a direct sum of representations $D':=\d'_1\t\d'_2\t ...\t\d'_k$ where $\d'_i$ are 
automorphic discrete series of $GL_{\frac{n}{k}}(\aa_\E)$. Each $D'$ is irreducible. If $U$ is the space of $D'$, then $\s^{-1}:U\to \s^{-1}U$ intertwines the representation $(D'^{\s^{-1}},U)$ with the representation by right translation in $\s^{-1}U$, that is  
${\d'_1}^{\s^{-1}}\t{\d'_2}^{\s^{-1}}\t ...\t{\d'_k}^{\s^{-1}}$. So $\s^{-1}:U\to \s^{-1} U$ intertwines $(D',U)$ with 
$({\d'_1}^{\s^{-1}}\t{\d'_2}^{\s^{-1}}\t ...\t{\d'_k}^{\s^{-1}})^\s$.

Now the restriction of $M(s_k,0)$ to $\s^{-1}U$, induces an operator $M:\s^{-1}U\to W$, where $W$ is the space of 
${\d'_2}^{\s^{-1}}\t{\d'_3}^{\s^{-1}}\t ...\t{\d'_k}^{\s^{-1}}\t{\d'_1}^{\s^{-1}}$, 
 which intertwines ${\d'_1}^{\s^{-1}}\t{\d'_2}^{\s^{-1}}\t ...\t{\d'_k}^{\s^{-1}}$ with ${\d'_2}^{\s^{-1}}\t{\d'_3}^{\s^{-1}}\t ...\t{\d'_k}^{\s^{-1}}\t{\d'_1}^{\s^{-1}}$. We have $W=U$ if and only if $\d'_2={\d'_1}^\s$, $\d'_3={\d'_2}^\s$, ...,$\d'_1={\d'_k}^\s$, if and only if 
 $D'$ is of type $\Delta\times \Delta^\s\times \Delta ^{\s^2}\times ...\times \Delta ^{\s^{k-1}}$, for $\D$ an automorphic discrete series of $GL_{\frac{n}{k}}(\aa_\E)$ and $\D^{\s^k}=\D$. In this case, $M$ intertwines ${\d'_1}^{\s^{-1}}\t{\d'_2}^{\s^{-1}}\t ...\t{\d'_k}^{\s^{-1}}$ with $D'$, hence $({\d'_1}^{\s^{-1}}\t{\d'_2}^{\s^{-1}}\t ...\t{\d'_k}^{\s^{-1}})^\s$ with ${D'}^\s$, and $M\circ \s^{-1}$ intertwines $(D',U)$ with $({D'}^\s,U)$.

The decomposition of  $\rho_{L'_k,t,\mu_\E}$ in direct sum of irreducible subrepresentations gives rise to a decomposition of its space in a direct sum of subspaces all stable by right translation and which are permuted by the operators $M(s,0)$ and $\s^{-1}$. 
So $\tr(M_{L'}^{GL_{n}(\aa_{\E})}(s,0)\ \circ\ \s^{-1}\ \circ \rho_{L',t,\mu_\E}(\phi))$ may be computed by taking the restriction of this operator to subrepresentations of type $\Delta\times \Delta^\s\times \Delta ^{\s^2}\times ...\times \Delta ^{\s^{k-1}}$ of $\rho_{L',t,\mu_\E}(\phi)$, i.e. those defined in a space $U$ such that $M(s,0)\circ \s^{-1}(U)=U$.

Then, computing explicitly the coefficients (see [Ba1], page 414) we come to the equality:

\begin{equation}\label {spectral}
\sum_{k|n}\frac{1}{k^2}
 \sum_{\d
 }
\tr(M_{L_k}^{GL_{n}(\aa_{\F})}(s_k,0)\ \circ \d^k(f))=
\end{equation}
$$l\, \sum_{k|n}\frac{1}{k^2}
 \sum_{\Delta
 }
\tr(M_{L'_k}^{GL_{n}(\aa_{\E})}(s_k,0)\ \circ \s^{-1}\ \circ \Delta\times \Delta^\s\times \Delta ^{\s^2}\times ...\times \Delta ^{\s^{k-1}}(\phi))
$$
where

- for $k|n$, $\d$ runs over the set of automorphic discrete series
of $GL_{\frac{n}{k}}(\aa_{\F})$  such that 

\ \ \ \ \ \ \ \ \ \ \ \ \ \ \ \ - $\d$ is $\mu'$-equivariant for some
character $\mu'$ of $\F^\t_\infty$ such that $\mu'^k=\mu$ and 

\ \ \ \ \ \ \ \ \ \ \ \ \ \ \ \ - the norm
of the imaginary part of its infinitesimal character is
$\frac{t}{k}$, 

\def\l{\lambda}

- for $k|n$, $\Delta$ runs over the set of automorphic discrete series of $GL_{\frac{n}{k}}(\aa_{\E})$ which are 

\ \ \ \ \ \ \ \ \ \ \ \ \ \ \ \ - $\s^{k}$-stable and

\ \ \ \ \ \ \ \ \ \ \ \ \ \ \ \ - $\mu_\E'$-equivariant for some character $\mu'_\E$ of $\E^\t_\infty$ such that $\mu_\E'^k=\mu_\E$, and 

\ \ \ \ \ \ \ \ \ \ \ \ \ \ \ \ - such that the norm
of the imaginary part of the infinitesimal character of $\Delta$ is
$\frac{t}{k}$,

- $\phi$ and $f$ are associated,

- $s_k$ is the Weyl element associated to the cyclic permutation $(k,k-1,...,1)$ of blocks.\\
\ \\
\ \\
{\bf Simplification.}

Set $s:=s_k$.  
For some representation $D:=\Delta\t\Delta^\s\times \Delta ^{\s^2}\times ...\times \Delta ^{\s^{k-1}}$ as in the formula, because 
$M(s,0)\circ \s^{-1}$ intertwines $(D,U)$ with $({D}^\s,U)$, and because $D$ is irreducible, $M(s,0)\circ \s^{-1}$ is, by Schur's lemma, a non-zero scalar multiple of the global intertwining operator $I_D$ obtained from local ones. We write $(M(s,0)\circ\s^{-1})_{|U}=\lambda_D I_D$, $\l_D\in \cc^*$. However, when $k=1$ (and so $M(s,0)$ is trivial), we have seen in section \ref{notglobal} that $\l_D=1$ (Proposition \ref{opintMW}).


\begin{lemme}\label{trivial}
If $D:=\Delta\t\Delta^\s\times \Delta ^{\s^2}\times ...\times \Delta ^{\s^{k-1}}$ and $D':=\s D=\Delta^\s\times \Delta ^{\s^2}\times ...\times \Delta ^{\s^{k-1}}\t\Delta$, then $\l_{D'}=\l_D$. 
\end{lemme}
\ \\
{\bf Proof.} The relation  $M(\tau^{-1},\tau\pi_E)\s^{-1}=\s^{-1} M(\tau^{-1},\tau\s \pi_E)$ from [AC] page 208 (we replaced $\s$ by $\s^{-1}$), specialized to $\tau=s^{-1}$, gives 
$$(M(s,0)\circ\s^{-1})_{|U}=(\s^{-1}\circ M(s,0))_{|U}.$$ Then 
$$(M(s,0)\circ\s^{-1})_{|U}=(\s^{-1}\circ M(s,0)\circ\s^{-1}\circ\s)_{|U},$$ 
hence 
$$\l_D I_D= \s^{-1}_{|U}(\l_{D'}I_{D'})\s_{|U}=\l_{D'}\ (\s^{-1}_{|U}I_{D'}\s_{|U}).$$ 
where $\s_{|U}:U\to \s U$. As $\s_{|U}$ intertwines $D^\s$ with $D'$, we have $\s^{-1}_{|U}I_{D'}\s_{|U} = I_{D^{\s}}$ by Lemma \ref{nio}. But $I_{D^\s}=I_D$, so the result follows.\qed


\ \\

Let $\pi_\F$ be an automorphic discrete series of $GL_n(\aa_{\F})$. 
Let $S$ be a finite set of places of $\F$ containing the infinite places, the finite places $v$ where, for $w|v$, $\E_w/\F_v$ is ramified, and the places where $\pi_\F$ is not spherical. 
Let $V$ be the set of places  of $\F$ complementary to $S$. 
Let $X_{\F,\pi_\F}$ be the set of representations $\pi'$ on the left of  Formula \ref{spectral} such that $\pi'_v$ has the same Shintani lift as $\pi_{\F,v}$ for all $v\in V$. 
By Proposition \ref{separating}, this set is finite, contains only automorphic discrete series, and is equal to $\{\chi\otimes \pi_\F,\ \chi\in X(\E/\F)\}$. We denote its cardinality $m(\pi_\F)$ or simply $m$.
Let $X_{\E,\pi_\F}$ be the set of representations $\Pi$ on the right such that $\Pi_v=\otimes_{w|v}\Pi_w$ is a Shintani lift of $\pi_{\F,v}$ for all $v\in V$ (i.e. $\Pi_w$ is a Shintani lift of $\pi_{\F,v}$ if $w|v$). 
{\it For every $\Pi\in X_{\E,\pi_\F}$, fix an isomorphism of $\Pi$ onto the restricted product ${\otimes}_v \Pi_v$}, $v$ place of $\F$, and denote, for every place $v$ of $\F$, by $I_{\Pi_v}$ the intertwining operator of $\Pi_v$ with $\Pi_v^\s$ induced by the action of $\s$ on $\Pi$, which is the \nio\ of $\Pi_v$ as shown at section \ref{notglobal}.

\begin{lemme} We have:
\begin{equation}\label{utile}
\sum_{\pi'\in X_{\F,\pi_\F}} \prod_{v\in S} \tr (\pi'_v(f_v))=
\end{equation}
$$l\, \sum_{k|n}\frac{1}{k^2}
 \sum_{\Pi\in X_{\E,\pi_\F}} \l_\Pi \prod_{v\in S}
\tr(I_{\Pi_v}\, \Pi_v(\phi_v))
$$
if $f_v\lra \phi_v\ \forall v\in S$.
\end{lemme}
\ \\
\ \\
{\bf Proof.} Let $v\in V$. Let $w_1,w_2,...,w_k$ the places of $\E$ dividing $v$ and, for $1\leq i\leq k$, identify $\E_{w_i}$ with $\E_{w_1}$ using $\s$-action. If $\phi_v=\otimes_{1\leq i\leq k}\phi_k$ is in the spherical Hecke algebra of $GL_n(\E)_v:=\prod_{w|v}GL_n(\E_w)$, then $f_v=b(\phi)$ -- where $\phi=\phi_1*\phi_2*...*\phi_k$ -- by definition, 
and if $\phi_i$ are all spherical, then $f_v$ is spherical ([AC], page 49). We  
 will apply Formula \ref{spectral} only with this kind of spherical functions at places $v\in V$ (and $w|v$, $v\in V$). The formula involves then only representations of $GL_n(\aa_{\F})$ which are spherical (and unitary) at all places $v\in V$ and representations of $GL_n(\aa_{\E})$ which are spherical (and unitary) at all places $w|v$ with $v\in V$. 
Each unitary spherical representation $\pi'_v$ of $GL_n(\F_v)$, for $v$ finite, admits a local Shintani lift $\Pi_v$ ({\bf Theorem B}, already proved in section \ref{sphericalunit}). For $v\in V$, if $\phi_i$ and $f_v$ are spherical, defined as before,
 we replace $\tr (\pi'_v(f_v))$ with $\tr (I_{\Pi_v}\Pi_v(\phi_v))$.
This equals
 $\tr (I_{\Pi_{w_1}}\circ\Pi_{w_1}(\phi))$. For $v\in V$ and $w|v$, $\tr (I_{\Pi_w}\circ\Pi_w(\phi))=\tr (\Pi_w(\phi))$ 
({\bf Remark}, section \ref{sphericalunit}) so the intertwining operators $I_{\Pi_w}$ do not play any role in the traces at these places and the standard method ([La], [JL], [Fl2]) applies to get the proposition.\qed
\ \\

Now, $X_{\F,\pi_\F}$ is equal to $\{\chi\otimes \pi_\F,\ \chi\in X(\E/\F)\}$ (Proposition \ref{separating}). If $\chi\in X(\E/\F)$ and $v$ is a place of $\F$, the character $\chi_v$ is trivial on the image of the norm map. Then, for 
$\pi'\in X_{\F,\pi_\F}$, the restriction of the character of $\pi'_v = \chi_v\otimes \pi_{\F,v}$ to the image of the norm map equals the restriction of the character of $\pi_{\F,v}$. 
If $f$ is as in Proposition \ref{intorb}, the regular orbital integral of $f_v$ is zero outside the image of the norm map. We have then: 
$$\tr (\pi'_{v}(f_v))=\tr (\pi_{\F,v}(f_v)).$$
As a consequence, the left side of the formula is equal to 
$$m\prod_{v\in S} \tr (\pi_{\F,v}(f_v)).$$ 
(``The image of norm map" means ``the set of elements of $G_{\F,v}$ which are conjugate to a norm".)

As one may choose functions $f_v$ such that this quantity does not vanish, we have that $X_{\E,\pi_\F}$ is not empty.
Let $D:=\Delta\times \Delta^\s\times \Delta^{\s^2}\times ...\times \Delta^{\s^{k-1}}\in X_{\E,\pi_\F}$. 
If $D'=\Delta'\times \Delta'^\s\times \Delta'^{\s^2}\times ...\times \Delta'^{\s^{t-1}}$ is another element of $X_{\E,\pi_\F}$, then $D'$ is isomorphic at almost every place (all the places outside $S$) with $D:=\Delta\times \Delta^\s\times \Delta^{\s^2}\times ...\times \Delta^{\s^{k-1}}$. By the strong multiplicity one theorem in the automorphic spectrum ([JS2]), which says that if two automorphic representations are isomorphic at almost every place they have the same cuspidal support, we have $\Delta'=\Delta^{\s^i}$ for some $i$, $1\leq i\leq r(\D)$. So  $X_{\E,\pi_\F}$ has $r(\D)$ elements, precisely the $D_i=\Delta^{\s^i}\t\D^{\s^{i+1}}\t \D^{\s^{i+2}}\t...\t\D^{\s^{i+k-1}}$ for $0\leq i\leq r(\D)-1$. By Lemma \ref{trivial}, $\l_{D_i}=\l_D$.  The right side of the equality \ref{utile} is then:
$$\l_D\, l\, \frac{r(\D)}{k^2}\prod_{v\in S}\tr(I_{v} D_v(\phi_v)),$$
where $I_v$ is the \nio\ of $D_v$, and we used Lemma \ref{nio} for the equality of traces associated to representations $D_i$, $0\leq i\leq r(\D)-1$.

We have proved: 

\begin{lemme}
There exist $k|n$, a square integrable representation $\D$ of $GL_{\frac{n}{k}}(\aa_{\E})$ and a complex number $\lambda_D$ such that, if we set  $D:=\Delta\times \Delta^\s\times \Delta^{\s^2}\times ...\times \Delta^{\s^{k-1}}$, then $D_v$  is the Shintani lift of $\pi_{\F,v}$ for $v\in V$ and, moreover,
\begin{equation}\label{simplified}
m\prod_{v\in S} \tr (\pi_{\F,v}(f_v))=\l_D\, l\, \frac{r(\D)}{k^2}\prod_{v\in S}\tr(I_{v} D_v(\phi_v)).
\end{equation}
if $\phi_v\lra f_v$ for all $v\in S$. 
\end{lemme}
\ \\

The following lemma will play some role in the proof.

\begin{lemme}\label{modulusone}
The operators $(M(s,0)\circ\s^{-1})_{|U}$ and $I_D$ are unitary. The modulus of the complex number $\l_D$ is $1$.
\end{lemme}
\ \\
{\bf Proof.} 
The subrepresentation ${\mathcal D}:=\oplus_{0\leq i\leq r(\D)-1} D_i$ of $\rho_{L',t,\mu_\E}$, $L'=L_k(\aa_\E)$, 
is stable by both $\s$ and $M(s,0)$. Let us show that $\s$ is a unitary operator. If we set 
$J_i=\Delta^{\s^i}\otimes\D^{\s^{i+1}}\otimes \D^{\s^{i+2}}\otimes ...\otimes \D^{\s^{i+k-1}}$ 
for $0\leq i\leq r(\D)-1$, then $D_i$ is induced from $J_i$ and ${\mathcal D}$ is induced from ${\mathcal J}:=\oplus_{0\leq i\leq r(\D)-1} J_i$. To show that $\s$ is an unitary operator of ${\mathcal D}$ it is enough to show that the action of $\s$ on ${\mathcal J}$ is unitary. The representations $J_i$ are in direct orthogonal sum, and $\s(J_i)=J_{i+1}$ for $0\leq i\leq r(\D)-1$ with the convention $J_{r(\D)}=J_0$. If $\o_i$ is the central character of $J_i$, then the scalar product on $J_i$ is defined by 
$$(f,h)_k=\int_{L_k(\E)Z(L_k(\aa_\E))\bc L_k(\aa_\E)}\overline{f(\bar{g})}h(\bar{g})d\bar{g},$$
restriction from the one of $L^2(L_k(\E)Z(L_k(\aa_\E))\bc L_k(\aa_\E);\omega_i)$.
Notice that $\s$ stabilizes the measure on $L_k(\E)Z(L_k(\aa_\E))\bc L_k(\aa_\E)$ because a finite order automorphism stabilizes a Haar measure. So it is easy to see that $\s$ is an isometry from $J_i$ to $J_{i+1}$.

The global intertwining operator $M(s,^.)$ is known to be unitary at $0$ ([MW2], IV.3.12). So the composed global operator $M(s,0)\circ\s^{-1}$ is also unitary.\\ 

Let us show that $I_D$ is unitary. Let $(\ ,\ )$ be the scalar product on $D$. Then $(\ ,\ )$ is stable not only by $D$ but also by $D^\s$. The representations $D$ and $D^\s$ are unitary and we identify their spaces with their duals using $(\ ,\ )$. If $I_D^*$ is the adjoint of $I_D$ for (\ ,\ ), then $I_D^*\circ I_D$ is an intertwining operator of $D$. Because $D$ is irreducible, $I_D^*\circ I_D=\l Id$, $\l\in\cc$  (by Schur's lemma). Because $(I_D^*\circ I_D(v), v)=(I_D(v),I_D(v))$, $\l$ is real positive. Then $\frac{1}{\sqrt{\l}} I_D$ is unitary. But $I_D^l=Id$, so $\l=1$.\\

The relation $(M(s,0)\circ \s^{-1})_{|U}=\lambda_D I_D$ is then of type $Y_1=\l_D Y_2$ with $Y_1$ and $Y_2$ unitary operators. So  $|\l_D|=1$.\qed
\ \\

Let us recall how, starting from this point, [AC] deals with the automorphic cuspidal case, and then adapt it to the case of residual representations in the next section. We start from Lemma \ref{simplified} with its notation.\\
\ \\
{\bf Three important steps.}

(1) The lemma implies that $\pi_\F$ is cuspidal if and only if $\D$ is cuspidal. Indeed, if $\pi_\F$ is cuspidal, it is generic at a finite place $v$ outside $S$. Then $D$ is generic at any place $w$ which divides $v$. Then so is $\D$ (easy from local classifications), and an automorphic discrete series with a generic component is cuspidal. 

Conversely, if $\D$ is cuspidal, then $\D$ is generic at a finite place $w$ outside $S$. Then $D$ is generic at $w$. So $\pi_\F$ is generic at the place $v$ which is divisible by $w$, which implies $\pi_\F$ is cuspidal.\\

(2) Recall that we are in the situation where we don't know yet the lift of representations other than spherical unitary representations (the other proofs were based on the lift of square integrable representations). 
In [AC] I Sec. 6, Arthur and Clozel use the simple twisted trace formula to prove first a rough version of the lift of a square integrable representation 
[AC], I  Th. 6.2.a,b. In the case of a local cyclic extension $E/F$, let us say that an irreducible representation $\pi$ of $GL_n(F)$ {\bf has a weak lift} if there exist a $\s$-stable irreducible representation $\Pi$ of $GL_n(E)$ such that $\pi$ and $\Pi$ verify the Shintani relation.
The proof (page 56) does not depend on $l$ being prime or not, and shows that a square integrable representation always has a weak lift, {\it which is a local component of a  global cuspidal representation of} $GL_n(\aa_\E)$. In particular, every square integrable representation of $GL_n(F)$  ($F$ local) has a weak lift, {\it which is a unitary irreducible generic representation}. Now every generic representation is an irreducible product of essentially square integrable representations, and every product of unitary generic representations is again an irreducible unitary generic representation. 

A consequence of Proposition \ref{liftprod} is then:

{\bf (G)} {\it In the local setting, every unitary generic representation $\pi$ of $GL_n(F)$ has a weak lift $\Pi$ which is a unitary generic representation of $GL_n(E)$}. (This is part of {\bf Theorem A} (c)).\\

(3) We go back to the relation \ref{simplified}. For now, the representation $\pi_\F$ is cuspidal or residual.
We know $r(\D)$ divides $k$ (because $\Delta\times \Delta^\s\times \Delta^{\s^2}\times ...\times \Delta^{\s^{k-1}}$ is $\s$-stable). Also, $\Delta\times \Delta^\s\times \Delta^{\s^2}\times ...\times \Delta^{\s^{k-1}}$ may be written as $\Theta^{\frac{k}{r(\D)}}$, with $\Theta=\Delta\times \Delta^\s\times \Delta^{\s^2}\times ...\times \Delta^{\s^{r(\D)-1}}$.\\ 
\ \\
{\it We assume, by induction, that {\bf Theorem E} is true for $GL_t(\F)$ and $GL_t(\E)$ for all $t<n$.}\\
\ \\
Then, if $r(\D)<k$,  $\Theta$ is the Shintani lift of some discrete series $\theta_\F$ of $GL_{\frac{nr(\D)}{k}}(\F)$. Then $\pi_\F$ and the automorphic representation $\theta_\F^{\frac{k}{r(\D)}}$ have the same local Shintani lift at every place outside $S$. By Proposition \ref{separating},  
 $r(\D)=k$, which leads to contradiction. 
So, {\it if we assume {\bf Theorem E} for $t<n$, we have, in the lemma, $r(\D)=k$ and so}

\begin{equation}\label{equa4}
m\prod_{v\in S} \tr (\pi_{\F,v}(f_v))=\l_D\, \frac{l}{k}\prod_{v\in S}\tr(I_{v} D_v(\phi_v)).
\end{equation}
\ \\

Following these three steps, one (i.e. Arthur and Clozel, Theorem 5.1 page 212) gets {\bf Theorem E} for cuspidal representations (actually, the step (1) implies that, if we assume {\bf Theorem E} for cuspidal representations for $t<n$, then we still get the relation \ref{equa4} if $\pi_\F$ is cuspidal). Indeed, the weak lift for all unitary generic representations ({\bf (G)}, Step (2)) implies (when $\pi_\F$ and $\D$ are cuspidal, cf. Step (1)) the equality $m=\l_D\, \frac{l}{k}$ from Equation \ref{equa4}. So $\l_D$ is real positive, and, as it has modulus $1$ (Lemma \ref{modulusone}), we have $\l_D=1$ and $mk=l$. Recall $r(\D)=k$ (Step 3).

As a consequence one (Arthur and Clozel, Proposition 6.6 page 58) gets the {\bf Theorem A}. Claims (a) and (b) are shown using that a local unitary cuspidal or square integrable representation may be realized as a local component of an automorphic cuspidal representation. Claim (c) is then a consequence of 
Proposition \ref{liftprod} and the fact
 that tempered representation are irreducible induced representations from square integrable representations.\\
\ \\
\subsection{Shintani lift for local unitary representations and global residual representations (Theorem D and E)}

The proofs of these lifts are interdependent. From the trace formula (more precisely Formula \ref{simplified}) we will get the global lifts of some particular global residual representations. Specialized to one place, that will imply the Shintani lifts of Speh representations. Then we will get the Shintani lifts of all unitary representations, by the local Proposition \ref{liftprod}. The local Shintani lifts of all unitary representations imply then the Shintani lifts of all global residual representations.\\

We adapt the preceding method to the case of residual representations. We do not have the result analogous to the step (2), because the local component of residual representations is not generic. So to get the general result we will prove first the Shintani lift for all unitary representations.

To show it, we would like to know that $m(\pi_\F) k=l$ and $\l_D=1$ also for a residual representation $\pi_\F$. But these equalities have been obtained from \ref{equa4} for cuspidal representations because we knew that local Shintani lifts of their local components exist ({\bf (G)}). We solve these two points below  in the following way: 

1. We show the relation $m(\pi_\F) k=l$ using the construction of residual representations from cuspidal ones and the fact that the relation has been proved for cuspidal ones. 

2. Then we will use, first, the relation \ref{equa4} in a particular case when we know $k=1$, and, so $\l_D=1$, to get the local Shintani lift for Speh  representations. This implies the lift of all unitary irreducible representations. The local Shintani lift of all  irreducible unitary representations then
implies, using \ref{simplified}, the global lift for all residual representations.\\

1. Let $\rho_\F$ be a cuspidal representation of $GL_n(\aa_\F)$ and set $\pi_\F:=MW(\rho_\F,q)$. Let $V$ be the set of finite  places of $\F$ where $\rho_\F$ is spherical and $\E/\F$ is unramified, and $S$ be the complementary set of places. Let $R:=\rho_\E\t\rho_\E^\s\t...\t\rho_\E^{\s^{r(\rho_\E)-1}}$, where $r(\rho_\E)=\frac{l}{m(\rho_\F)}$, be the Shintani lift of $\rho_\F$.
 For every $v\in V$, $R_v$ is the Shintani lift (of spherical representations) of $\rho_{\F,v}$. Set $\Pi:=MW(\rho_\E,q)$. Set $P:= \Pi\t\Pi^\s\t...\t\Pi^{\s^{r(\rho_\E)-1}}$. For every $v\in V$, $P_v$ is the Shintani lift (of spherical representations) of $\pi_{\F,v}$ (this follows from the classification of unitary spherical representation and {\bf Theorem B}). So, if $D=\d\t\d^\s\t...\t\d^{\s^{k-1}}$ is the representation appearing in the equality \ref{equa4}:

$$m(\pi_\F)(\prod_{v\in S} \tr \pi_{\F,v}(f_v))=\l_D\, \frac{l}{k}\prod_{v\in S}\tr(I_{v} D_v(\phi_v)),$$
then $D=P$. The equality between cuspidal supports of $D$ and $P$ implies $k=r(\rho_\E)$ (check the central exponents). By unicity of the cuspidal support, again, one has $m(\pi_\F)=m(\rho_\F)$ (the characters stabilizing $\pi_\F$ up to isomorphism are the characters stabilizing $\rho_\F$ up to isomorphism). So the equation reads now:

\begin{equation}\label{equa5}
\prod_{v\in S} \tr \pi_{\F,v}(f_v)=\l_D\, \prod_{v\in S}\tr(I_{v} P_v(\phi_v)).
\end{equation}
\ \\

2. Let $\d$ be a square integrable representation of $GL_n(F)$, $F=\r$ or $F=\cc$ or $F$ a $p$-adic field. Assume $\F$ is chosen such that, for some place $v_0$, $\F_{v_0}\simeq F$ (this is always possible). Let us choose the cuspidal representation $\rho_\F$ with local component isomorphic to $\d$ at the place $v_0$ and  cuspidal at some other place $v_1$. 
The proof of the existence of such a representation is the same as in [AC] lemma 6.5, page 54. Let $S$ be the set consisting of the infinite places and all the finite places where $\rho_\F$ is not spherical.
Let $\E$ be a cyclic extension of $\F$ of degree $l$ such that $\E_{v_0}\simeq E$ and $\E$
splits at all the places in $S\bc\{v_0\}$ (for the existence of such an extension, Theorem 1.2.2 in [AC], quoting [AT]). Then $m(\rho_\F)=l$, because the lift $D$ of $\pi_\F$ is cuspidal at places dividing $v_1$ and so $k=1$. Then $\l_D$ is $1$ in equation \ref{equa5}, which becomes:

$$\prod_{v\in S} \tr \pi_{\F,v}(f_v)=\prod_{v\in S}\tr(I_{v} P_v(\phi_v)).$$

The Shintani lift at $v\in S\bc\{v_0\}$ is trivial since $\E_v$ splits, so the Shintani relation at the place $v_0$ follows. 
But the local component of $\pi_\F$ at the place $v_0$ is, by construction, $u(\d,q)$. 
And the local component of $P$ is $\Pi_{v_0}$. 
Now, $\Pi_{v_0}=u(\rho_{\E,v_0},q)$, where $\rho_{\E,v_0}$ is generic, equal to the local base change of $\d$ as in {\bf Theorem A} (b). 
This proves {\bf Theorem D} (a).\\ 

Let us prove {\bf Theorem D} (b).
Using Proposition \ref{liftprod} and the local classification of unitary representations ([Ta1], [Ta2]), the Shintani lift for representations $u(\d,q)$ ({\bf Theorem D} (a)) implies the Shintani lift of all unitary representations. 
Let now $u$ be a $\s$-stable unitary representation of $GL_n(E)$. 
Then one may write $u=\prod_i u_i$ where $u_i$ are representations of type $u(D,q)$ or $\pi(u(D,q);\alpha)$ with $D$ square integrable and 
$0<\alpha<\frac{1}{2}$. 
By unicity of the terms of the product, $u^\s=\prod_i u_i^\s$. As $u$ is $\s$-stable, it is not hard to see that $u=\prod_j U_j$, where $U_j$ are representations of type $u(D,q)\t u(D,q)^\s\t...\t u(D,q)^{\s^{r(D)-1}}$ or $\pi(u(D,q),\alpha)\t \pi(u(D,q),\alpha)^\s\t...\t \pi(u(D,q),\alpha)^{\s^{r(D)-1}}$. 
As each of these representations are Shintani lift of unitary representations of $GL_n(F)$ (by {\bf Theorem D} (a)). 
By Proposition \ref{liftprod}, $u$ is the Shintani lift of some unitary representation of $GL_n(F)$.\\

Let us show {\bf Theorem E} for automorphic discrete series.
The local lift of unitary representations, re-injected in Equation \ref{equa4} implies the global Shintani lift for discrete series. Indeed, we now know the local Shintani lifts for all local components. This shows (a).\\
\ \\

We prove {\bf Theorem E} (b). Let $\Pi=\pi_\E\tt  \pi_\E^\s \tt...\tt  \pi_\E^{\s^{r-1}}$ be a representation of $GL_n(\aa_{\E})$, where $\pi_\E$ is an automorphic discrete series of $GL_{\frac{n}{r}}(\aa_{\E})$ such that $r(\pi_\E)=r$. Starting with Formula \ref{traceformula}, we let $V$ be the places $v$ where $\Pi_v$ is spherical and $\E_w/\F_v$ is unramified for $w|v$. Let $S$ be the complementary set of places of $\F$. Let $X_{\E,\Pi}$ be the set of representations $\Pi'$ of $GL_n(\E)$ appearing on the right side of Formula \ref{traceformula} such that $\Pi'_v\simeq \Pi_v$ for all $v\in V$. Let $X_{\F,\Pi}$ be the set of representations $\pi'$ appearing on the left side of Formula \ref{traceformula} such that $\Pi_v$ is a Shintani lift of $\pi'_v$ for all $v\in V$. Then by the same simplification argument we have already used, we get the relation:

\begin{equation}\label{utile2}
\sum_{\pi'\in X_{\F,\Pi}} \prod_{v\in S} \tr (\pi'_v(f_v))=
\end{equation}
$$l\, \sum_{k|n}\frac{1}{k^2}
 \sum_{\Pi'\in X_{\E,\Pi}} \l_{\Pi'} \prod_{v\in S}
\tr(I_{\Pi'_v}\, \Pi'_v(\phi_v))
$$
if $f_v\lra \phi_v$.

By the same arguments we have already explained, $X_{\E,\Pi}$ has $r=r(\pi_\E)$ elements, all having the same contribution. Moreover, $r(\pi_\E)=k$ by the assumption before equation \ref{equa4}. We get:
\begin{equation}\label{utile3}
\sum_{\pi'\in X_{\F,\Pi}} \prod_{v\in S} \tr (\pi'_v(f_v))=
\end{equation}
$$\frac{l}{r}
 \l_{\Pi} \prod_{v\in S}
\tr(I_{\Pi_v}\, \Pi_v(\phi_v))
$$
if $f_v\lra \phi_v$. As the right-hand side is not identically zero (locally, the trace distribution of an irreducible representation is not identically zero so we may always find a function $\phi_v$ such that $\tr(I_{\Pi_v}\, \Pi_v(\phi_v))\neq 0)$), the left hand side is not identically zero, and $X_{\F,\Pi}$ is not empty. Then it contains at least one element $\pi$. By the same arguments as before, we prove that $\Pi_v$ is the Shintani lift of $\pi_v$. This proves {\bf Theorem E} (b) for discrete series.\\
\ \\
{\bf Remark.} Shintani lift is well understood for generic, spherical, elliptic and unitary representations (as well as their twists with characters). 
Proposition \ref{liftprod} allows one to lift products of such representations, under some conditions of irreducibility. In his recent paper [Ta4], Tadi\'c gives a simple criterion to know when a product of twists of unitary representations is irreducible (one implication has been proved by different methods in [MW1]).

\newpage
\section{Appendix: multiplicity one irreducible subquotients}

Intertwining operators between representations induce intertwining operators between isomorphic multiplicity one irreducible subquotients. This is what we will study in this section.
We start with a few definitions in Algebra used in the text.

We need some definitions concerning representations of a general group $G$. We are only interested here in ``complex" representations and we omit the term complex in what follows.

Let $(\Pi,V)$ be a representation of $G$. If $U,W$ are subspaces of $V$ which are stable by $\Pi$, we denote $\Pi_U$ the subrepresentation of $\Pi$ in $U$ and $\Pi_{U/W}$ the quotient representation of $\Pi_U$ in $U/W$ induced by $\Pi$. Then $\Pi_{U/W}$ is said to be a {\bf subquotient} of $\Pi$. 
We say that $\Pi$ has {\bf finite length} if there exists a finite sequence $0=V_0\subset V_1\subset V_2\subset ...\subset V_k=V$ of subspaces of $V$ stable by $\Pi$ such that the subquotient representation $\Pi_{V_i/V_{i-1}}$ is irreducible for all $1\leq i\leq k$. 
It is known, [Bou] I.1.4.7, that if $\Pi$ is of finite length, then for any other such sequence $0=V'_0\subset V'_1\subset V'_2\subset ...\subset V'_{k'}=V$ such that $\Pi_{V'_i/V'_{i-1}}$ is irreducible for $1\leq i\leq k'$ there is a bijection $\s$ of $\{1,2,...,k\}$ onto $\{1,2,...,k'\}$ such that $\Pi_{V_i/V_{i-1}}$ is isomorphic to $\Pi_{V'_{\s(i)}/V'_{\s(i)-1}}$ for $1\leq i\leq k$. Then $k=k'$ is called the {\bf length} of $\Pi$. If $\tau$ is an irreducible representation of $G$, its {\bf multiplicity} in $\Pi$ is the number $r$ of indexes $i\in\{1,2,...k\}$ such that $\tau$ is isomorphic to $\Pi_{V_i/V_{i-1}}$. We say that $\tau$ {\bf appears} in $\Pi$ if its multiplicity in $\Pi$ is $\geq 1$. Writing $Irr(G)$ for the set of isomorphism classes of irreducible representations of $G$, we get from $\Pi$ a map  $JH(\Pi):Irr(G)\to\n$, where $JH(\Pi)(\tau)$ is the multiplicity of $\tau$ in $\Pi$ for $\tau\in Irr(G)$. If $0\to \Pi\to \Pi'\to \Pi''\to 0$ is an exact sequence, then we have the equality $JH(\Pi')=JH(\Pi) + JH(\Pi'')$.\\

Let $\Pi$ be a representation of finite length of $G$ and 
$\tau$ be an irreducible subquotient of $\Pi$ of multiplicity one. Let $X$ be the set of pairs $(U',W')$ of stable subspaces of $V$ such that $W'\subset U'$ and $\Pi_{U'/W'}\simeq \tau$. Choose $(U,W)\in X$ with $U$ maximal.

\begin{prop}\label{subquot}
{\rm (a)} Let $(U',W')\in X$. Then $U'\subset U$, $W'\subset W$ and $W'=W\cap U'$. 
Moreover the inclusion $U'\subset U$ induces an isomorphism $\Pi_{U'/W'}\simeq \Pi_{U/W}$.

{\rm (b)} If $(U',W')\in X$ with $U'$ maximal, then $(U',W')=(U,W)$.

{\rm (c)} If $f$ is an automorphism of $\Pi$, then $f(U)=U$ and $f(W)=W$.
\end{prop}

{\bf Proof.} (a) Let $(U',W')\in X$. The obvious map $U/W\to (U+U')/(W+U')$ is surjective; as $\Pi_{U/W}\simeq \tau$, $\Pi_{(U+U')/(W+U')}$ is either zero or isomorphic with $\tau$. In the second case, $U'\subset U$ by the maximality of $U$. Moreover, $W+U'=W$, so $U'\subset W$; but then $\tau$ appears in $\Pi_W$ because it appears in $\Pi_{U'}$, and this contradicts the multiplicity one assumption. So we are in fact in the first case, where $U+U'=W+U'$, and, in particular, $U\subset W+U'$. Consider the natural exact sequence: 
$$0\longrightarrow W/(W\cap W')\longrightarrow (W+U')/W'\longrightarrow (W+U')/(W+W')\longrightarrow 0.$$
Now $\tau$ appears in $\Pi_{(W+U')/W'}$ because $\Pi_{U'/W'}\simeq \tau$, so $\tau$ appears either in $\Pi_{W/(W\cap W')}$ or in $\Pi_{(W+U')/(W+W')}$. 
But $\tau$ cannot appear in $\Pi_{W/(W\cap W')}$ because $\Pi_U$ contains $\tau$ with multiplicity one, 
and $\Pi_{U/W}\simeq \tau$. Then $\tau$ appears in $\Pi_{(W+U')/(W+W')}$. 
As we have seen that $U\subset W+U'$, we get $U=W+U'$ by maximality of $U$, and in particular $U'\subset U$. 
Then the quotient map $U/W\to U/(W+W')$ is an isomorphism, so $W'\subset W$. 
But the quotient map $\Pi_{U'/W'} \to \Pi_{U/W}$ is an isomorphism, because both representations are isomorphic to $\tau$, and consequently $W'=U'\cap W$. 
The final assertion in (a) is an immediate consequence, and so is (b). 



(c) Let $f$ be an automorphism of $\Pi$. Then $f$  induces an isomorphism of representations $\Pi_{U/W}\simeq \Pi_{f(U)/f(W)}$, so $(f(U),f(W))\in X$. By (b), $f(U)\subset U$ and $f(W)\subset W$. Applying the same argument to the isomorphism $f^{-1}$ we get (c).\qed
\ \\
\ \\
{\it {\bf Remarks.} {\rm (1)} In the situation of (c), $f$ induces an automorphism $\bar{f}$ of the representation $\Pi_{U/W}$. Choosing a $G$-isomorphism $\phi:\tau\simeq \Pi_{U/W}$ yields an automorphism $f_\tau$ of $\tau$ such that $\bar{f}\circ\phi =\phi\circ f_\tau$; if Schur's lemma is valid for $\tau$, the automorphism $f_\tau$ does not depend on the choice of $\phi$.

{\rm (2)} One may take from the beginning $U$ minimal instead of $U$ maximal and the Proposition and its proof hold similarly. We will not use that in our paper.}\\

Assume now an intertwining operator $f:\Pi\to \Pi'$ is given, where $\Pi'$ is another representation of $G$. If one fixes an isomorphism $\tau\simeq \Pi_{U/W}$, one gets by composition with $f$ an intertwining operator $f':\tau\to \Pi_{f(U)/f(W)}$. Then $f'$ is an isomorphism if and only if the image of $f$ contains a subquotient isomorphic to $\tau$. We then call the intertwining operator $f'$ the restriction of $f$ to $\tau$. 
 If one starts with  $\tau \simeq \Pi_{U'/W'}$, with $(U',W')\in X$ not maximal, then the canonical morphisms 
 $\Pi_{U'/W'}\to \Pi_{U/W}$ and $\Pi'_{f(U')/f(W')}\to \Pi'_{f(U)/f(W)}$, which turn out to be here isomorphisms, induce an isomorphism $\Pi_{U'/W'}\to \Pi'_{f(U')/f(W')}$ which is the same as the one induced by $f$.\\
\ \\
{\bf Group with automorphism.}
In the main text the group $G$ is endowed with an automorphism $\s$; we now
explore that context. For a representation $(\Pi,V)$ of $G$ we get another one $(\Pi^\s,V)$, where $\Pi^\s(g):=\Pi(\s(g))$ for all $g\in G$. Assume $\Pi$ is of finite length and $\tau$ is an irreducible subquotient of $\Pi$. A subspace of $V$ is stable by 
$\Pi$ if and only if it is stable by $\Pi^\s$. Moreover, if $U,W$ are stable subspaces of $V$ such that $W\subset U$, we have $(\Pi_U)^\s=(\Pi^\s)_U$ and $(\Pi_{U/W})^\s=(\Pi^\s)_{U/W}$.
So :

(a) the set of pairs $(U,W)$ where $U,W$ are subspaces of $(\Pi,V)$ such that $\tau\simeq \Pi_{U/W}$ equals the set of pairs $(U',W')$ where $U',W'$ are subspaces of $(\Pi^\s,V)$ such that $\tau^\s \simeq \Pi^\s_{U'/W'}$; if $v$ is the space of $\tau$, an isomorphism $\phi:v\to U/W$ which intertwines $\tau$ and $\Pi_{U/W}$ intertwines $\tau^\s$ and $\Pi^\s_{U/W}$.

Assume $\Pi\simeq \Pi^\s$ and fix an isomorphism $f:\Pi\to\Pi^\s$. Then :

(b) $(U,W)\mapsto (f(U),f(W))$ induces a bijection from the set of pairs $(U,W)$
where $U,W$ are subspaces of $(\Pi,V)$ such that $\tau\simeq \Pi_{U/W}$ onto the set of pairs $(U',W')$ where $U',W'$ are subspaces of $(\Pi^\s,V)$ such that $\tau^\s \simeq \Pi^\s_{U'/W'}$.\\

Assume now $\tau$ has multiplicity one in $(\Pi,V)$ and $\tau$ is isomorphic to $\tau^\s$. 
Let $(U,W)$ be the unique maximal pair such that $\tau\simeq \Pi_{U/W}$ as in the Proposition \ref{subquot}. 
Fix an isomorphism $\phi:\tau\simeq \Pi_{U/W}$. By (a), (b), and the Proposition \ref{subquot} (b) we have $f(U)=U$, $f(W)=W$, and 
$f$ induces the quotient isomorphism $\bar{f}:\Pi_{U/W}\to \Pi^\s_{U/W}$. 
We get the intertwining operator $\phi^{-1}\bar{f}\phi:\tau\to \tau^\s$, which does not depend on $\phi$ if Schur's lemma applies to $\tau$. 
We say then that {\it $\phi^{-1}\bar{f}\phi$ is the $\s$-operator on $\tau$ obtained from $f$ by the multiplicity one property} (it is an operator on the space of $\tau$ which intertwines $\tau$ with $\tau^\s$).
Note that this situation is more general than the claim (c) of the Proposition, where $f$ was an automorphism of $\Pi$. 

\newpage

\section{Bibliography}

[AC] J.Arthur, L.Clozel, {\it Simple Algebras, Base Change, and the Advanced Theory of the Trace Formula},  
Ann. of Math. Studies, Princeton Univ. Press 120, (1989).

[AT] E. Artin, J. Tate, {\it Class Field Theory}, Benjamin, New York, 1967.

[Ba1] A. I. Badulescu, Global Jacquet-Langlands correspondence, multiplicity one and classification
of automorphic representations, Invent. Math., 172(2):383 - 438 (2008). With an appendix by Neven Grbac.

[Ba2] A. I. Badulescu, Jacquet-Langlands et unitarisabilit\'e,  J. Inst. Math. Jussieu 6 (2007), no. 3, 349 - 379.

[Bar] E. M. Baruch, A proof of Kirillov's conjecture. Ann. of Math. (2) 158 (2003), no. 1, 207 - 252.

[Be] J. N. Bernstein, P-invariant distributions on $GL(N)$ and the classification of unitary
representations of $GL(N)$, (non-Archimedean case), in {\it Lie groups and representations II}, Lecture
Notes in Mathematics 1041, Springer-Verlag, 1983.

[BJ], A. Borel, H. Jacquet, Automorphic forms and automorphic representations.  Proc. Sympos. Pure Math., XXXIII, 
{\it Automorphic forms, representations and L-functions}, Part 1, pp. 189-207, Amer. Math. Soc., Providence, R.I., 1979. 

[BZ] J. N. Bernstein, A. V. Zelevinsky, Induced representations of reductive $p$-adic groups I, Ann. Sci. \'Ec. Norm. Sup. (4) 10 (1977), no. 4, 441 - 472.

[Bou] N. Bourbaki, {\it Algebra I (Chapters 1-3)}, Springer 1989. 

[BW] A. Borel, N. Wallach, {\it Continuous cohomology, discrete subgroups, and representations of reductive groups}, Annals of Mathematics Studies, 94. Princeton University Press, Princeton, N.J.; University of Tokyo Press, Tokyo, 1980. xvii+388 pp. 

[Cl] Clozel, L. Th\'eor\`eme d'Atiyah-Bott pour les vari\'et\'es $p$-adiques et caract\`eres des groupes r\'eductifs, in  
{\it Harmonic analysis on Lie groups and symmetric spaces} (Kleebach, 1983). M\'em. Soc. Math. France (N.S.) No. 15 (1984), 39 - 64. 

[DS] S. DeBacker and P. J. Sally, Jr., {\it Admissible invariant distributions on reductive $p$-adic groups},
University Lecture Series, 16. American Mathematical Society,  Providence, RI, 1999.

[Fl1] D. Flath, Decomposition of representations into tensor products, Proc. Sympos. Pure Math., XXXIII, {\it Automorphic forms, representations and L-functions}, Part 1, pp. 179-185, Amer. Math. Soc., Providence, R.I., 1979. 

[Fl2] D. Flath, A comparison for the automorphic representations of $GL(3)$ and its twisted
forms, Pacific J. Math. 97 (1981), 373 - 402.

[H-C]  Harish-Chandra,
Admissible invariant distributions on reductive $p$-adic groups,
{\it Lie theories and their applications}, (Proc. Ann. Sem. Canad. Math.  Congr., Queen's Univ., Kingston, Ont., 1977), pp. 281-347, 
Queen's  Papers in Pure Appl. Math., No. 48, Queen's Univ., Kingston, Ont., 1978.

[He1] G. Henniart, Induction automorphe globale pour les corps de nombres,  Bull. Soc. Math. France 140 (2012), no. 1, 1 - 17.

[He2] G. Henniart, Induction automorphe pour $GL(n,\cc)$. J. Funct. Anal. 258 (2010), no. 9, 3082 - 3096. 

[HL] G. Henniart, B. Lemaire, Changement de base et induction automorphe pour ${\rm GL}_n$ en caractéristique non nulle,  M\'em. Soc. Math. Fr. (N.S.) No. 124 (2011).

[JL] H. Jacquet, R. P. Langlands, {\it Automorphic forms on $GL(2)$}, Lecture Notes in Math. 114,
Springer-Verlag (1970).

[Ja] H. Jacquet, On the residual spectrum of ${\rm GL}(n)$, in {\it Lie group representations II} (College Park, Md., 1982/1983), 185 - 208, 
Lecture Notes in Math., 1041, Springer, Berlin, 1984.

[JP-SS] H. Jacquet, I. Piatetski-Shapiro, J.A.Shalika, Rankin-Selberg convolutions, Amer. J. Math. 105 (1983), no. 2, 367-464. 

[JS1] H. Jacquet, J. Shalika, The Whittaker models of induced representations, Pacific J. Math. 109 (1983), no. 1, 107-120.

[JS2] H.Jacquet, J.A.Shalika, On Euler products and the classification
of automorphic forms II, Amer. J. Math.  103 (1981), no. 4, 777 - 815.

[Ko] B. Kostant, On Whittaker vectors and representation theory, Invent. Math. 48 (1978), no. 2, 101 - 184.


[La] R.P.Langlands, {\it Base change for $GL(2)$},  Ann. of Math. Studies,
Princeton Univ. Press 96, (1980).

[La2] Langlands On the notion of automorphic representation,
Proc. Sympos. Pure Math., XXXIII, {\it Automorphic forms, representations and L-functions}, Part 1, pp. 208-215, Amer. Math. Soc., Providence, R.I., 1979. 

[MW1] C. Moeglin, J.-L. Waldspurger, Le spectre r\'esiduel de $GL(n)$, Ann. Sci. \'Ec. Norm. Sup., t 22, (1989), 605 - 674.

[MW2] C. Moeglin, J.-L. Waldspurger, {\it Spectral decomposition and Eisenstein series}, Cambridge Tracts in Mathematics, 113. 
Cambridge University Press, Cambridge, 1995.

[P-S] I. Piatetski-Shapiro, Multiplicity one Theorems, Proc. Sympos. Pure Math., XXXIII, {\it Automorphic forms, representations and L-functions}, Part 1, pp. 209-212, Amer. Math. Soc., Providence, R.I., 1979.

[Re] D. Renard, {\it Repr\'esentations des groupes r\'eductifs p-adiques}, Cours sp\'ecialis\'e SMF, {\bf 17} (2011).

[Ro] F. Rodier,
Whittaker models for admissible representations of reductive p-adic split groups, Proc. Symp. Pure Math., XXVI,
{\it Harmonic analysis on homogeneous spaces},  pp. 425-430. Amer. Math. Soc., Providence, R.I., 1973.

[Sh] J. A. Shalika, The multiplicity one Theorem for $GL_n$, Ann. of Math. 100 (1974), 171 - 193.

[S1] F. Shahidi, On certain L-functions, Amer. J. Math. 103 (1981), no. 2, 297 - 355.

[S2] F. Shahidi, Fourier transforms of intertwining operators and Plancherel measures for GL(n), Amer. J. Math. 106 (1984), no. 1, 67 - 111. 

[Ta1] M. Tadi\'c, Classification of unitary representations in irreducible representations of general linear group (non-Archimedean case), Ann. Sci. \'Ec.  Norm. Sup. (4) 19 (1986), no. 3, 335 - 382.

[Ta2] M. Tadi\'c, $GL(n,\cc)^\wedge$ and $GL(n,\r)^\wedge$, {\it Automorphic forms and L-functions II: local aspects}, 285-313, Contemporary Mathematics, vol. 489, Am. Math. Soc., Providence, RI (2009).

[Ta3] M. Tadi\'c, Induced representations of ${\rm GL}(n,A)$ for $p$-adic division algebras $A$, J. Reine Angew. Math. 405 (1990), 48-77.

[Ta4] M. Tadi\'c, Irreducibility criterion for representations induced by essentially unitary ones (case of non-archimedean GL(n,A)), 
Glasnik Matematicki, vol. 49, no. 1, (2014). 

[Vo] D. A. Vogan Jr., Gel'fand-Kirillov dimension for Harish-Chandra modules, Invent. Math. 48 (1978), no. 1, 75-98.

[Vo2] D. A. Vogan Jr., The unitary dual of $GL(n)$ over an Archimedean field, Invent. Math. 83 (1986), 449 - 505.

[Wa] N. Wallach, {\it Real reductive groups II}, Pure and Applied Mathematics, 132-II. Academic Press, Inc., Boston, MA (1992), xiv+454 pp.

[Ze]  A.Zelevinsky, Induced representations of reductive $p$-adic groups II,  Ann. Sci. \'Ec. Norm. Sup. 13 (1980), 165 - 210.\\
\ \\
\ \\
\ \\
\ \\
A. I. Badulescu, Universit\'e de Montpellier 2, Institut Montpelli\'erain Alexander Grothendieck,\\
Case Courrier 051, Place Eug\`ene Bataillon, 34095 Montpellier Cedex, France\\
E-mail : ibadules@univ-montp2.fr\\
\ \\
G. Henniart, Universit\'e de Paris-Sud, Laboratoire de Math\'ematiques,\\
Orsay cedex F-91405 France, CNRS, Orsay cedex F-91405 France\\
E-mail : Guy.Henniart@math.u-psud.fr

\end{document}